\documentclass[onecolumn]{autart_muga}  
\usepackage{lineno,hyperref}
\modulolinenumbers[5]

\usepackage{graphicx}
\usepackage{amssymb,amsmath}
\usepackage{lipsum}
\usepackage{mathtools}
\usepackage{cuted}
\usepackage{breqn}
 \usepackage{enumerate}
\newtheorem{theorem}{Theorem}
\newtheorem{corollary}{Corollary}

\newtheorem{assumption}{Assumption}
\newtheorem{remark}{Remark}
\usepackage{lineno,hyperref}
\usepackage{graphicx}
\graphicspath{ {images/}}
\begin{document}
\begin{frontmatter}
\title{A Geometric PID Control Framework for \\ Mechanical Systems}

\author[mymainaddress]{D.H.S. Maithripala},
\author[mysecondaryaddress]{T.W.U. Madhushani},
\author[mytertiaryaddress]{J.M. Berg}
\address[mymainaddress]{Dept. of Mechanical Engineering, Faculty of Engineering, University of Peradeniya, KY 20400, Sri Lanka.\\  smaithri@pdn.ac.lk}
\address[mysecondaryaddress]{Postgraduate and Research Unit, Sri Lanka Technological Campus, CO 10500, Sri Lanka.\\  udarim@sltc.lk}
\address[mytertiaryaddress]{Dept. of Mechanical Engineering, Texas Tech University, TX 79409, USA.\\ jordan.berg@ttu.edu}

\begin{abstract}
The inherent coordinate independent nature of a geometric controller makes it an ideal candidate for achieving the best global stability properties that one can expect for a given physical system.
These lectures demonstrate the development of  such a PID control framework for mechanical systems. The starting point of this effort is the observation that mechanical systems are intrinsically double integrator systems. For a double integrator system on a linear space a PD controller, ensures the globally exponential stability of the origin, by giving the closed loop system the structure of a spring mass damper system. The incorporation of Integral action makes the closed loop stability properties robust to bounded unmodelled constant input disturbances and bounded parametric uncertainty. We generalize this linear PID controller to mechanical systems that have a non-Euclidean configuration space. Specifically we start by presenting the development of the geometric PID controller for fully actuated mechanical systems
and then extend it to a class of under actuated interconnected mechanical systems of practical significance by introducing the notion of \emph{feedback regularization}. We show that feedback regularization is the mechanical system equivalent to partial feedback linearization. 
We apply these results for trajectory tracking for several systems of interest in the field of robotics. First, we demonstrate the robust almost-global stability properties of the geometric PID controller developed for fully actuated mechanical systems using simulations and experiments for a multi-rotor-aerial-vehicle. The extension to the class of under actuated interconnected systems allow one to ensure the semi-almost-global locally exponential tracking of  the geometric center of a spherical robot on an inclined plane of unknown angle of inclination. The results are demonstrated using simulations for a hoop rolling on an inclined plane and then for a sphere rolling on an inclined plane. The final extension that we present here is that of geometric PID control for holonomically or non-holonomically constrained mechanical systems on Lie groups. The results are demonstrated by ensuring the robust almost global locally exponential tracking of a nontrivial spherical pendulum.
\end{abstract}
\end{frontmatter}


\section{Introduction}
It is often sufficient to represent the state of a dynamic system as a vector in Euclidean space, and to treat the change of the state with respect to time as a vector in the same Euclidean space. Many results in control theory have been developed in this framework. However there are also cases of practical importance where the system naturally evolves on a more general manifold, and in these cases a generalized control framework is desirable. For example, the set of possible orientations of a rigid body do not constitute a vector space, nor is this set topologically equivalent to any Euclidean space, nor is the angular velocity an element from the same space as the orientation. A well-developed analytical framework is available to describe dynamics of systems on smooth manifolds, however tools for control on these spaces are less advanced. The objective of this work has been to extend the useful concept of ``integral action'' from Euclidean space to a more general setting in the most natural way.

Smooth manifolds are described by a collection of overlapping ``charts,'' each of which is diffeomorphic, through a typically nonlinear set of coordinate functions, to an open subset of $\mathbb{R}^n$ for some constant $n$. Roughly speaking, there is an increasingly general hierarchy of control methods for systems on manifolds. First, the dynamics may be linearized about a single point on the manifold, and linear control techniques employed. Second,  a controller may be designed for the nonlinear equations within a single coordinate chart. Finally, a controller may be designed for the manifold itself, free of any specific choice of coordinate chart. This last approach is sometimes referred to as ``geometric'' or ``intrinsic.'' With each increase in generality comes the potential for improved stability and performance over a larger region of the state space. These factors become important when the operating region cannot be restricted in advance, as for vehicles undergoing aggressive maneuvers or for systems that have an arbitrary initial configuration. Furthermore, it has been shown that if the underlying configuration space is not diffeomorphic to $\mathbb{R}^n$, there exists no continuous state feedback that globally asymptotically stabilizes a given configuration \cite{Bhat}. The best that can be achieved using smooth feedback on a non-Euclidean space is almost-global stability. Such global topological constraints are inherently absent from any design that considers only a single coordinate chart. A controller that stabilizes an entire chart may seem to ensure global stability when in fact it does not. Thus a geometric approach gives a more accurate characterization of global stability.

An individual or a set of interacting rigid bodies moving under the influence of external forces is called a \textit{mechanical system}. Formally, such a system is characterized by its configuration space, the kinetic energy, the constraints, and the external forces. The behavior of these systems is governed by the Euler-Lagrange equations which are a generalization of Newton's equations. These equations are intrinsic. Meaning they do not depend on the choice of coordinates used to specify the configuration. However a particular choice of coordinates are necessary for the explicit expression of them.
Newton's equations for a free particle moving under the influence of a control force results in a double integrator system.
The proportional + integral + derivative, or PID, controller is the simplest robust controller that one can design for regulating a double integrator system. It ensures the globally exponential convergence of the tracking error for ramp  references in the presence of constant input disturbances and bounded parametric uncertainty. Thus for a free particle the PID controller can be considered to be the ideal choice for ensuring robust tracking. The Euler-Lagrange equation describing the motion of a more complicated mechanical system comprising interconnected rigid bodies is expressed using coordinates. In general these equations do not take the form of a double integrator system in the conventional sense. The reason for this apparent absence of the double integrator structure is a consequence of the non-Euclidean nature of the configuration space of rigid body motion and the choice of coordinates one has to make in order to represent the Euler-Lagrange equations.

Typically the configuration space of a rigid body system takes the structure of a smooth non-Euclidean manifold. The kinetic energy of the system allows one to define an inner product on each of the tangent spaces to the manifold in a smooth manner. This assignment allows one to define a metric on the manifold. A manifold equipped such a metric structure is referred to as a Riemannian manifold \cite{Marsden,Bullo} The evolution of the system under the influence of the forces acting on the system describes a twice differentiable curve on the configuration manifold called the configuration trajectory. The tangent to this curve is an intrinsic quantity referred to as a tangent vector and is commonly referred to as the generalized velocity of the system. The space of all possible such velocities at a given configuration is the tangent space to the manifold at that configuration. The generalized forces acting on the system are found to be objects that belong to the dual space of the tangent space called the cotangent space. The ordinary derivative of a tangent vector taken with respect to a particular choice of coordinate system on the manifold, will in general, fail to be a tangent vector. Therefore the ordinary coordinate differentiation, of tangent vectors, does not correspond to an intrinsic operation. An intrinsic differentiation operation that ensures that the derivative of a tangent vector is also a tangent vector is called a covariant derivative operator. It is well known that there exists a unique such derivative operator called the Levi-Civita connection that is compatible with the Lie derivative on the manifold and expresses the rate of change of the kinetic energy of the system as the ``product'' between the force and the velocity. This unique covariant derivative turns out to be the correct intrinsic notion of acceleration that allows the expression of a mechanical system in a form that is equivalent to a double integrator.  Thus it may turn out to be possible to extend linear system concepts used for regulation of double integrator systems to mechanical systems in an intrinsic way. The work presented here is motivated by this observation.

The \emph{first major contribution} that is reported here is the extension of the notion of PID control to the class of fully actuated mechanical systems proposed in \cite{MaithripalaAutomatica}. The concept of proportional + derivative, or PD, control was first extended to the geometric setting in \cite{Koditschek}, and subsequently further developed in \cite{BulloTracking,MaithripalaLieGroup,MaithripalaSIAM,BanavarPD}. The central object in PD control is the tracking error, and the insight of the geometric extension is to give the tracking error dynamics the form of a {mechanical system}, or in other words the structure of a double integrator system on a Riemannian manifold. Geometric PD control design follows from deriving control action from artificial ``potential energy'' terms that create ``energy minima'' where the tracking error is zero. While the ideas of geometric PD control can be applied to general Riemannian manifolds, the expressions are much more compact for systems evolving on Lie groups. More recently, we have augmented geometric PD control by incorporating integral action to obtain true geometric PID control \cite{MaithripalaAutomatica}. The key insight in this extension is the observation made in the previous paragraph that mechanical systems are essentially double integrator systems on a Riemannian manifold. Implementing an integral error term based on this insight we  \cite{MaithripalaAutomatica} obtain robust configuration tracking of a fully actuated mechanical system on a Lie group. The integral action provides almost-global, locally exponential convergence of the tracking error to zero in the presence of bounded parametric uncertainty and bounded constant disturbance forces. In \cite{QuadICIIS2015} we experimentally demonstrate the excellent robust global stability properties of this controller for attitude tracking of a Multi Rotor Aerial Vehicle (MRAV). The result demonstrates the capability of recovering from very large deviations and, to date, is the best reported experimentally verified global stability property found in the literature.

The result in \cite{MaithripalaAutomatica} only applies to  configuration tracking of fully-actuated systems. The broader class of systems that we consider in \cite{RollingHoopACC2017,RollingBallAutomatica} correspond to two or more mechanical systems coupled through shared control forces. Considered independently, each subsystem is fully actuated, but since inputs are shared by two or more configuration variables, the combined system is underactuated. The subsystems may also be coupled through quadratic velocity terms arising from the Reimannian structure corresponding to the overall system kinetic energy. The \emph{second major contribution} reported here is to extend geometric PID control to output tracking of this class of underactuated mechanical systems \cite{RollingHoopACC2017}. 
Because the system is underactuated, it will not be possible to achieve configuration tracking of all configuration variables. Thus we have designated one of the mechanical systems as the {output system}, and collectively referred to the others as the {actuator system}. The control objective was to ensure that a output that depends only on the configuration of the output system tracks a desired output. We have assumed that the output is relative degree two with respect to the controls of the system and that the zero dynamics of the system has a stable but not necessarily  asymptotically stable equilibrium.
The {output error system} is then the system describing the discrepancy between the output of the system and the desired output reference trajectory. The goal of the underactuated trajectory tracking problem is to achieve asymptotic convergence of the output error to zero, while ensuring that the actuator system remains stable. 

While the class of coupled systems described above may seem overly narrow, it arises naturally when multi-body systems interact with each other. In such situations the inputs are shared owing to the fact that the interaction forces and moments must be equal and opposite. 
Due to these interactions the expression of the dynamics of 
the individual subsystems may fail to correspond to a mechanical system or in other words to an intrinsic double integrator structure on the configuration space of the individual subsystem. However we observe that feedback control may be used to make them look like interconnected mechanical systems. In some sense, this is a geometric interpretation of partial feedback linearization (PFL). PFL is a powerful technique of nonlinear control, in which state feedback is used to cancel all nonlinearities in the system input-output response \cite{Isidori}. In contrast rather than cancel terms we introduce terms that are quadratic in the velocities so that the individual systems will take the form of mechanical systems or equivalently look like intrinsic double integrator systems. Unlike the PFL procedure, our feedback terms (the terms we add) are independent of coordinate systems, and therefore can be used to provide the best possible global stability results. {The major contribution of the work we present in \cite{RollingHoopACC2017} is to use coordinate-independent feedback to inject quadratic velocity terms that correspond to the Levi-Civita connection for the system kinetic energy and thereby provide each of the subsystems with the structure of a mechanical system}. Since the objective is not a linear system, but rather a mechanical system, we refer to this process as {feedback regularization}. 
Under the assumptions of the paper \cite{RollingHoopACC2017}, every system to which feedback regularization can be applied can subsequently be controlled to track a desired output trajectory using geometric PID. {Therefore the other major contribution of \cite{RollingHoopACC2017} is to achieve general output trajectory tracking for the class of underactuated mechanical systems described above.}

Spherical robots provide an ideal test bed to demonstrate these techniques. The coupling between the spherical body and the actuation mechanism is through reaction forces and moments. If no other forces were present, Newton's laws applied to the error system augmented by the coupled system would produce equations of motion suitable for the geometric PID tracking controller. However, the forces corresponding to the interaction between the subsystems destroy that structure. Thus we first use feedback regularization to recover the Riemannian structure of the subsystems, and subsequently apply geometric PID to the regularized system. These ideas are first demonstrated on a hoop rolling without slip on an inclined surface in \cite{RollingHoopACC2017} and subsequently extended to a spherical robot rolling without slip on an inclined surface of unknown angle of inclination in \cite{RollingBallAutomatica}. We show that the controller is capable of ensuring that the center of mass of the sphere tracks a twice differentiable reference trajectory almost-semi-globally with local exponential convergence in the presence of bounded constant unmodelled disturbances and parameter uncertainty while ensuring the stability of the actuation mechanism. This result significantly extends the hitherto known work on tracking for spherical robots and is the \emph{third major contribution} that is reported here.

The extension of the PID controller developed for Lie groups to Riemannian manifolds have two major difficulties. The first one is the absence of a convenient intrinsic notion of configuration error while the next major difficulty is the construction of an intrinsic velocity error. The first problem may be overcome by considering distance functions on the configuration manifold or by embedding the configuration space in a higher dimensional Euclidean space \cite{Chaturvedi,BanavarPD} while the second problem may be overcome by considering the parallel transport map that is induced by the Levi-Civita connection \cite{BulloTracking,Chaturvedi,BanavarPD}. In here we report an alternative approach that exploits the observation that the natural configuration space of nontrivial mechanical systems is a finite product of the Euclidean motion group $SE(3)$ and hence is a Lie group. It is the presence of holonomic constraints that may cause a reduction of the configuration space to that of a configuration space that does not have the structure of a Lie group.
The most famous example is the nontrivial spherical pendulum. Being a rigid body rotating about a pivot point the natural configuration of the system is the group of rigid rotations, $SO(3)$. The holonomic constraint that prevents the pendulum from spinning about its axis reduces the configuration space of the pendulum to the sphere $\mathbb{S}^2$ that does not have the structure of a Lie group. However when one takes the constraint forces into account the system can be considered as a mechanical system on the Lie group $SO(3)$ instead of as an unconstrained system on $\mathbb{S}^2$. Motivated by this observation, the final and the \emph{fourth major contribution} that we report in here is that of extending the geometric control framework to holonomically or non-holonomically constrained systems on Lie groups \cite{UdariCDC2017}.  
We demonstrate the ideas by developing and simulating the  controller for a spherical pendulum. We show in \cite{UdariCDC2017} that the controller can be used to ensure the almost-global and locally exponentially stability of the upright configuration of a Segway type vehicle. 
This extension may have applications in formation control of mobile robots and un-manned aerial vehicles \cite{MaithripalaASME,UnderactuatedVehicle,FormationControlICIIS2013} and remains as a topic to be investigated.

This report is organized as follows. In Section-\ref{Secn:CovDerivative} we briefly review some of the basic notions of differential geometry necessary for expressing mechanical systems as double integrators. The framework we present is valid for both holonomically and non-holonomically constrained systems. The initial presentation is valid for singular-Riemannian manifolds. In Section-\ref{Secn:LieGroups}
we specialize it to the case where the configuration space is a Lie Group. On Lie Groups the expressions do not require coordinates due to the trivial nature of its tangent bundle. In Section-\ref{Secn:ExamplesMechSystems}, as examples, we will explicitly write down the expression of a general unconstrained mechanical system on the two commonly encountered Lie Groups: the circle $\mathbb{S}$ and the group of 3-dimensional Euclidean motion $SE(3)$. The major contributions of our work is reported in Section-\ref{Secn:PIDonLieGroups}. Section-\ref{Secn:PIDfullyActuated} presents the development of the intrinsic PID controller proposed in \cite{MaithripalaAutomatica} for configuration tracking of fully actuated mechanical systems and Section-\ref{Secn:NLPID} presents the intrinsic PID controller we propose in \cite{RollingHoopACC2017,RollingBallAutomatica} for a class of under actuated interconnected mechanical systems.  The Section-\ref{Secn:PID4Constrained} presents our final major contribution which is the extension of the geometric PID to holonomically or non-holonomically constrained mechanical systems on Lie groups.
These results are then demonstrated in Section-\ref{Secn:ExamplesPID_LieGroups}. Specifically in \ref{Secn:AttitudeStabilization} we present simulation and experimental verification results for attitude tracking of a MRAV.
From a bench mark point of view we also show in Section- \ref{Secn:IPC} that the controller can be used to ensure the almost-semi-global stabilization of the vertically upright equilibrium of an Inverted Pendulum on a Cart (IPC) moving on an inclined plane while ensuring that the cart velocities remain bounded. The result also has practical significance since an IPC is an approximation of a Segway moving with negligible yaw motion. The simulation verification of the PID controller that was developed in Section-\ref{Secn:NLPID} for interconnected mechanical systems is presented; first for a hoop rolling on an inclined plane in Section-\ref{Secn:RollingHoop} and then for a sphere rolling on an inclined plane in Section-\ref{Secn:RollingSphere}. Finally in Section-\ref{Secn:SphericalPendulum} we will demonstrate the geometric PID for constrained mechanical systems for robust almost global tracking of a spherical pendulum.

\section{Intrinsic Mechanical Systems}\label{Secn:CovDerivative}
A mechanical system is defined by a configuration space, the kinetic energy, the generalized forces, and the holonomic and non-holonomic constraints acting on the system \cite{Marsden,Bullo,Koditschek}. The generalized velocities and the generalized forces turn out to be elements of the tangent bundle and the cotangent bundle respectively while the kinetic energy induces the structure of a Riemannian manifold on the configuration space. 
In this context the Euler-Lagrange equation can be shown to take the form of a nonlinear double integrator. When the system is subjected to holonomic or nonholonomic constraints the forces that ensure the constraints do no work. This allows a reduction in the dimension of the double integrator system. In the following we briefly review some of the basic notions of differential geometry that allows us make precise these statements. 

Specifically, the configuration space $G$ of a mechanical system will be assumed to be a smooth $n$-dimensional manifold. The tangent vector to a curve through $g$ is denoted by $v_g=\dot{g}$.
Denote by $T_gG$ the space of all such tangent vectors to $G$ at $g\in G$. The collection of all such tangent spaces to $G$ is referred to as the tangent bundle, $TG$.  The space dual to $T_gG$ is denoted by $T_g^*G$ and is referred to as the cotangent space. The collection of all such cotangent spaces to $G$ is referred to as the cotangent bundle, $T^*G$.
It is customary to denote by $\langle\cdot,\cdot\rangle :
T_g^*G \times T_gG \mapsto \mathbb{R}$ the action of a covector $\gamma_g\in T_g^*G$ on a vector $v_g\in T_gG$ explicitly by $\langle \gamma_g,v_g\rangle$. 
The evolution of the mechanical system over time describes a smooth curve, $g(t)$, on the configuration space and the tangent to this curve $v_g=\dot{g}\in T_gG$ is the generalized velocity of the system while the generalized force acting on the system denoted by $\gamma_g$ is an elements of the cotangent space $T^*_gG$. Then $\langle \gamma_g,v_g\rangle$ is the rate of change of energy of the mechanical system. 
 
A singular Riemannian metric on $G$ assigns in a smooth fashion a degenerate inner product, 
$\langle\langle \cdot,\cdot \rangle \rangle : T_g G\times T_g G \mapsto \mathbb{R}$, on each of the tangent spaces $T_g G$ and gives $G$ the structure of a singular Riemannian manifold.  For mechanical systems the kinetic energy of the system defines such a metric so that the kinetic energy of the mechanical system is given by $\mathrm{KE}=\langle\langle \dot{g},\dot{g} \rangle \rangle/2$.
The singular Riemannian metric allows one to define a map $\mathbb{I}_g\::\:T_gG\mapsto T^*_gG$ by the relationship $\langle\mathbb{I}_gv_g,u_g\rangle\triangleq \langle\langle v_g,u_g\rangle\rangle$ for all $v_g,u_g\in T_gG$. The smooth tensor field $\mathbb{I}\::\:TG\mapsto T^*G$ that is point wise defined  above is usually referred to as the \textit{inertia tensor}.  For  Riemannian metrics the above defined map $\mathbb{I}_g$ is an isomorphism and in this case one can uniquely identify a vector with a given covector in an intrinsic fashion.

For any vectorfields $X,Y,Z\in TG$ the derivative of $\langle\langle X,Y\rangle\rangle$ along solutions of $Z$ is denoted by $\mathcal{L}_Z\langle\langle X,Y\rangle\rangle$.
For a Riemannian or singular Riemannian metric one can show that there exists a unique 1-form field, $\mathbb{I}\nabla_XY\in T^*G$ that satisfies the following properties \cite{Stocia}:
\begin{align}
\mathcal{L}_Z\langle\langle X,Y\rangle\rangle&=\langle \mathbb{I}\nabla_ZX,Y\rangle+\langle \mathbb{I}\nabla_ZY,X\rangle,\label{eq:Metric}\\
\mathbb{I}\nabla_XY-\mathbb{I}\nabla_YX&=\mathbb{I}[X,Y].\label{eq:Symmetric}
\end{align}
For a given $X,Y\in TG$ this 1-form field is explicitly given by the Koszul formula:
\begin{align*}
\langle\mathbb{I}\nabla_XY,Z\rangle&=\frac{1}{2}\left(\mathcal{L}_X\langle\langle Y, Z\rangle\rangle+\mathcal{L}_Y\langle\langle Z, X\rangle\rangle-\mathcal{L}_Z\langle\langle X,Y\rangle\rangle
-\langle\langle X,[Y, Z]\rangle\rangle+\langle\langle Y, [Z,X]\rangle\rangle+\langle\langle Z, [X,Y]\rangle\rangle\right).
\end{align*}
\noindent
It can be shown that this allows one to define a covariant derivative called the \textit{lower derivative} \cite{Stocia}, $\mathbb{I}\nabla : TG\times TG \mapsto T^*G$, that takes values in $T^*G$. Condition (\ref{eq:Metric}) states that the lower derivative is metric and (\ref{eq:Symmetric}) states that the  lower derivative is symmetric or torsion free. From a mechanical system point of view what turns out to be more crucial is the property of metricity given by (\ref{eq:Metric}).
If the metric is Riemannian then $\mathbb{I}$ is an isomorphism and then $\nabla_X Y\triangleq \mathbb{I}^{-1}(\mathbb{I}\nabla_XY)$ defines a unique \emph{covariant derivative} or \emph{connection},  called the Levi-Civita connection. 

\subsection{Unconstrained Mechanical Systems}
Using the above notations one can write down Euler-Lagrange equations for an \textit{unconstrained mechanical system} in the intrinsic fashion 
\begin{align}
\mathbb{I}\nabla_{\dot{g}}\dot{g}&=\gamma_g,\label{eq:MechanicalSystem}
\end{align} 
where $\gamma_g\in T^*_gG$ is the generalized force acting on the system. This expression is valid for both Riemannian metrics as well as singular Riemannian metrics. Since $\mathbb{I}$ is an isomorphism for  Riemannian (non-singular) metrics the mechanical system can be alternatively written as $\nabla_{\dot{g}}\dot{g}=\mathbb{I}^{-1}\gamma_g=\Gamma_g$.  Where as for singular Riemannian metrics the corresponding mechanical system can not be written in this fashion. In the Riemannian case $\nabla_{\dot{g}}\dot{g}$ has the notion of intrinsic acceleration.
When the configuration space is $\mathbb{R}^n$ and the kinetic energy metric is constant equation (\ref{eq:MechanicalSystem}) reduces to 
$m\ddot{x}=f$. Thus equation (\ref{eq:MechanicalSystem}) can be considered as the nonlinear equivalent of the linear double integrator.

\subsection{Constrained Mechanical Systems}
If the system is subjected to holonomic or non-holonomic constraints the Lagrange-d'Alembert principle states that the forces that ensure the constraints do no work. 
This observation allows a reduction in dimension of the representation of the double integrator system (\ref{eq:MechanicalSystem}).
At each $g\in G$ consider the collection of subspaces $\mathcal{D}_g\subset T_gG$. If
the assignment of these subspaces is smooth over $g$ the collection of all such subsaces, denoted by $\mathcal{D}\subset TG$, is called a smooth distribution. The holonomic and non-holonomic constraints of a mechanical system are typically of the form $\dot{g}(t)\in{\mathcal{D}_{g(t)}}$ for all $t$ for some such smooth $\mathcal{D}$ of dimension $m$. 
Let $\mathcal{D}^*_c\subset T^*G$ be the \emph{annihilating co-distribution} of $\mathcal{D}$ or simply the \emph{constraint co-distribution}. That is let 
\[
\mathcal{D}^*_c\triangleq\{\gamma\in T^*G\:\mid \langle \gamma ,v\rangle=0, \:\forall\: v\in \mathcal{D}\}.
\] 
The dimension of $\mathcal{D}^*_c$ is $(n-m)$. We will assume that $\mathbb{I}$ restricted to $\mathcal{D}$ is one-to-one and onto its image $\mathcal{D}^*\triangleq \mathbb{I}(\mathcal{D})$. Thus $\mathrm{dim}(\mathcal{D})=\mathrm{dim}(\mathcal{D}^*)=m$.
Since if $\dot{g}\in \mathcal{D}$ then $\mathbb{I}\dot{g} \notin \mathcal{D}^*_c$ it follows that $T^*G=\mathcal{D}^*\oplus \mathcal{D}^*_c$.

\begin{figure}[h!]
	\centering
	\begin{tabular}{c}
		\includegraphics[width=0.235\textwidth]{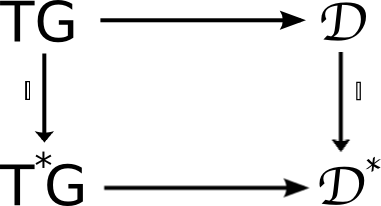}
	\end{tabular}
	\caption{Projection.\label{Fig:Projections}}
\end{figure}

Let $P_{\mathcal{D}^*}:T^*G\mapsto \mathcal{D}^*$ be the smooth projection onto $\mathcal{D}^*$. Thus $P_{\mathcal{D}^*_c}\triangleq (I-P_{\mathcal{D}^*})$ is the complementary projection that projects onto $\mathcal{D}^*_c$.  Let us summarize some of the consequences of this decomposition of $T^*G$.
\begin{enumerate}[(a)]
\item $P_{\mathcal{D}^*}(\mathbb{I}X)=\mathbb{I}X$ if and only if $X\in \mathcal{D}$.
\item For any $X\in \mathcal{D}$ the inner product $\langle \langle X,Y\rangle\rangle=0$ if and only if $Y \notin \mathcal{D}$.
\end{enumerate}

Let $\gamma (g)\in T^*_gG$ be the resultant generalized external force acting on the system and  $\gamma_\lambda(g)\in T^*_gG$ be the
generalized constraint force that maintains the constraint $\dot{g}(t)\in{\mathcal{D}_{g(t)}}$. For brevity of notation we will omit the $g$ dependence of the above forces. The Lagrange-d'Alembert principle says that the generalized constraint forces do no work and hence that $\gamma_\lambda\in \mathcal{D}^*_c$. Thus the 
Lagrange-d'Alembert principle implies that $P_{\mathcal{D}^*}(\gamma_{\lambda})\equiv 0$. The constraint that $\dot{g}\in{\mathcal{D}_{g(t)}}$ implies that $\mathbb{I}\dot{g}\notin{\mathcal{D}^*_c}$ and hence that $P_{\mathcal{D}^*_c}(\mathbb{I}\dot{g})\equiv 0$ and $P_{\mathcal{D}^*}(\mathbb{I}\dot{g})=\mathbb{I}\dot{g}$.

For a mechanical system corresponding to a singular-Riemannian metric, the intrinsic Newtons equations (\ref{eq:MechanicalSystem}) can be split along the constraint co-distribution $\mathcal{D}^*_{c}$ and its complement $\mathcal{D}^*$. Using this decomposition it can be shown that
the constrained equations of motion, of a mechanical system corresponding to a semi-Riemannian metric, have the intrinsic representation
\begin{align*}
P_{\mathcal{D}^*}(\mathbb{I}\nabla_{\dot{g}} {\dot{g}}) &= P_{\mathcal{D}^*}(\gamma), 
\\
P_{\mathcal{D}^*_c}(\mathbb{I}{\dot{g}}) &=0.
\end{align*}
Also
{
\begin{align}
\gamma_\lambda\!\!&=\!\!P_{\mathcal{D}^*_c}(\mathbb{I}\nabla_{\dot{g}}\dot{g}) - P_{\mathcal{D}^*_c}(\gamma)=-(\nabla_{\dot{g}}P_{\mathcal{D}^*_c})(\mathbb{I}\dot{g}) - P_{\mathcal{D}^*_c}(\gamma).\label{eq:ConstrainetForce}
\end{align}
}
Where the lower covariant derivative of $P_{\mathcal{D}^*_c}$ along $X\in TG$ is defined by
\begin{align*}
(\nabla_XP_{\mathcal{D}^*_c})(\mathbb{I}Y)\triangleq \mathcal{L}_X\left(P_{\mathcal{D}^*_c}(\mathbb{I}Y)\right)-P_{\mathcal{D}^*_c}(\mathbb{I}\nabla_XY).
\end{align*}
Basically what expression (\ref{eq:ConstrainetForce}) implies is that the constraint force should be equal and opposite to the projection of the external forces along the constraint co-distribution plus the `gyroscopic' forces along the constraint co-distribution.  
Combining the above three equations we have that the constrained mechanical system on the Lie group $G$ takes the form
\begin{align}
\mathbb{I}\nabla_{\dot{g}} {\dot{g}} &= -(\nabla_{\dot{g}}P_{\mathcal{D}^*_c})(\mathbb{I}\dot{g})+P_{\mathcal{D}^*}(\gamma). \label{eq:ConstrainedMech}
\end{align}

\subsection{Mechanical Systems on Lie Groups}\label{Secn:LieGroups}
Mechanical systems consists of interconnected rigid bodies. Thus the configuration space of a generic mechanical system is a finite product of $SE(3)$.  
Therefore, without loss of generality, mechanical systems can be considered as constrained or unconstrained mechanical systems with a configuration space equipped with the structure of a Lie group.  Due to the trivial nature of the tangent bundle to a Lie group we will see that the Levi-Civita connection and hence the representation of the mechanical system takes a particularly a simple form that does not require coordinates on the Lie Group. Below we will briefly review these notions.

Denote by $e$ the identity element of the Lie group $G$. Let $\mathcal{G}\triangleq T_eG$ be the Lie algebra of $G$ and $\mathcal{G}^*$ be its dual space. For $\zeta \in \mathcal{G}$ the notation $g\cdot \zeta$ will be used to represent the left translation of $\zeta \in T_eG$ to 
$g\cdot \zeta \in T_gG$ using the derivative of the left-multiplication map on $G$ while 
$ \zeta\cdot g$ will be used to represent the right translation of $\zeta \in T_eG$ to 
$ \zeta \cdot g\in T_gG$ using the derivative of the right-multiplication map on $G$. Conversely the left translation of a velocity vector $v_{g}\in T_gG$ to $T_eG$ denoted by $\nu^L_g\triangleq g^{-1}\cdot v_{g}$ will be referred to as the 
\textit{left-velocity} while the right translation to $T_eG$ given by $\nu^R_g\triangleq v_{g}\cdot g^{-1}$ will be referred to as the \textit{right-velocity} corresponding to $v_g$.  
Since the left and right velocities are globally defined $(g,v_{g})\in TG$ can be identified with $(g,\nu^L_g)\in G\times \mathcal{G}$ or $(g,\nu^R_g)\in G\times \mathcal{G}$. 
Being equal to the rate of change of momentum, forces are intrinsically co-vectors and hence elements of the cotangent space $T^*_gG$.  Let $\gamma_g\in T^*_gG$. The left and right translations allow one to pull-back the force $\gamma_g$ acting at $g$ to  $\mathcal{G}^*$ using either the left translation map or the right translation map. 
Denote by ${\gamma}^L_g$ and ${\gamma}^R_g$ the pull-back of $\gamma_g$ using the left and right translations  respectively. These are explicitly given by  $\langle{\gamma}^L_g,\zeta\rangle\triangleq \langle{\gamma}_g,g\cdot \zeta\rangle$ and $\langle{\gamma}^R_g,\zeta\rangle\triangleq \langle{\gamma}_g,\zeta \cdot g\rangle$. Thus we see that $(g,\alpha_{g})\in T^*G$ can be identified with $(g,{\gamma}^L_g)\in G\times \mathcal{G}^*$ or $(g,{\gamma}^R_g)\in G\times \mathcal{G}^*$. 
Thus the tangent bundle, $TG$,  and the cotangent bundle, $T^*G$, are trivial. That is $TG\equiv G\times \mathcal{G}$ and $T^*G\equiv G\times \mathcal{G}^*$. Furthermore we see that there are two distinct ways of trivialization.

Let $\langle\langle \cdot,\cdot \rangle \rangle_G : TG\times TG\mapsto \mathbb{R} $ be a Riemannian metric on $G$.  The metric along with left translations on $G$ induces the isomorphism $\mathbb{I}^L_g : \mathcal{G} \to \mathcal{G}^*$ by the relationship 
$\langle\mathbb{I}^L_g \zeta,\eta\rangle\triangleq \langle\langle g\cdot \zeta,g\cdot \eta\rangle \rangle_G$ while the right translations on $G$ induces the isomorphism $\mathbb{I}^R_g : \mathcal{G} \to \mathcal{G}^*$ by the relationship 
$\langle\mathbb{I}^R_g \zeta,\eta\rangle\triangleq \langle\langle \zeta\cdot g,\eta\cdot g\rangle \rangle_G$ for each $g\in G$. These induced isomorphisms are called left and right inertia tensors of the system respectively. In general these maps are a function of $g$. If $\mathbb{I}^L_g$ (resp. $\mathbb{I}^R_g$)  is a constant then the metric is said to be left-invariant (resp. right-invariant).  

The two different velocity representations of the system arising due to the use of left or right translations on the group will give two different explicit representations for the Levi-Civita connection.
For a given vector fields $X(g)$ on $G$ the left-trivialisation allows us to define a map
$\xi^L:G\to \mathcal{G}$ by $\xi^L_g\triangleq g^{-1}\cdot X(g)$. Similarly the right-trivialisation allows us to define a map
$\xi^R:G\to \mathcal{G}$ using $\xi^R_g\triangleq X(g)\cdot g^{-1}$. Let $\eta^L_g\triangleq  g^{-1}\cdot Y(g)$ and $\eta^R_g=Y(g)\cdot g^{-1}$. Then one can define 
$ \nabla^L_{\xi^L_g} {\eta^L_g}\triangleq  g^{-1}\cdot(\nabla_{X(g)}{Y(g)})$ and $\nabla^R_{\xi^R_g} {\eta^R_g}\triangleq (\nabla_{X(g)}{Y(g)})\cdot g^{-1}$. For a metric that does not possess any invariance properties these representations of the connection are in general functions of $g$. 
For notational convenience we will drop the $g$ dependence as well as the superscript in our subsequent representations where it will be understood that, in explicit computations, one will have to use $L$ or $R$ when considering left velocities $\zeta^L=g^{-1}\cdot\dot{g}$ or right velocities $\zeta^R=\dot{g}\cdot g^{-1}$ respectively. 
 
We conclude this section by pointing out that if the Riemannian metric satisfies certain invariant properties then the Levi-Civita connection can be explicitly expressed without the use of co-ordinates on $G$ or $\mathcal{G}$. In this case the Levi-Civita connection can be expressed using the adjoint action $\mathrm{Ad} : G\times \mathcal{G} \to \mathcal{G}$ that is defined by  $\mathrm{Ad}_g\eta =g\cdot \eta \cdot g^{-1}$ for all $g\in G$ and $\eta\in \mathcal{G}$. The associated derivative map at the identity $\mathrm{ad} : \mathcal{G}\times \mathcal{G} \to \mathcal{G}$ will be given by
$\mathrm{ad}_\zeta \eta =[\zeta,\eta]$ where $[\cdot,\cdot]$ is the Lie bracket on $\mathcal{G}$ while the dual of this map will be denoted by $\mathrm{ad}^*$.

Specifically it can be shown that, for left-invariant metrics,
\[
\nabla_\xi \eta \triangleq d{\eta}(\xi)+\frac{1}{2}\left(\pm\,\mathrm{ad}_\xi  \eta - \mathbb{I}^{-1}\left(\mathrm{ad}^*_{\xi}{\mathbb{I}\eta} +\mathrm{ad}^*_{\eta}{\mathbb{I}\xi}  \right)\right),
\]
and for right-invariant metrics,
\[
\nabla_\xi \eta \triangleq d{\eta}(\xi)+\frac{1}{2}\left(\pm\,\mathrm{ad}_\xi  \eta + \mathbb{I}^{-1}\left(\mathrm{ad}^*_{\xi}{\mathbb{I}\eta} +\mathrm{ad}^*_{\eta}{\mathbb{I}\xi}  \right)\right),
\]
where the $+$ results in $\nabla^L_{\xi^L} {\eta^L}$ and the $-$ results in $\nabla^R_{\xi^R} {\eta^R}$. 
Here 
$
d\eta(\xi)\triangleq \left.\frac{d}{dt}\right |_{t=0}\eta(g\exp{\xi t}),
$
for left-velocities while
$
d\eta(\xi)\triangleq \left.\frac{d}{dt}\right |_{t=0}\eta((\exp{\xi t})g),
$
for right-velocities.
A metric that is both left and right invariant is said to be a bi-invariant metric. It can be shown that
for a bi-invariant metric $\left(\mathrm{ad}^*_{\xi}{\mathbb{I}\eta} +\mathrm{ad}^*_{\eta}{\mathbb{I}\xi}  \right)\equiv 0$ and hence it follows that the bi-invariant connection is given by
\[
\nabla_\xi \eta \triangleq d{\eta}(\xi)\pm\frac{1}{2}\mathrm{ad}_\xi  \eta,
\]
where once again the $+$ results in $\nabla^L_{\xi^L} {\eta^L}$ and the $-$ results in $\nabla^R_{\xi^R} {\eta^R}$.

Thus when the configuration space is a Lie group $G$ the  equation (\ref{eq:MechanicalSystem}) that describes the motion of an unconstrained system takes the form
\begin{align}
\nabla_\zeta \zeta &= \gamma, \label{eq:UnconstrainedLieGroup}
\end{align}
while 
the  equation (\ref{eq:ConstrainedLieGroup}) that describes the motion of a constrained system takes the form
\begin{align}
\mathbb{I}\nabla_{\zeta} {\zeta} &= -(\nabla_{\zeta}P_{\mathcal{D}^*_c})(\mathbb{I}\zeta)+P_{\mathcal{D}^*}(\gamma). \label{eq:ConstrainedLieGroup}
\end{align}
When one chooses left-velocities to represent a system we will call it a \textit{left-velocity representation} while if right-velocities are used we will call it a \textit{right-velocity representation}.

\subsection{Examples}\label{Secn:ExamplesMechSystems}
For illustration purposes, in the following we consider in detail two specific classes of mechanical systems of practical significance: a mechanical system on a circle $\mathbb{S}$ and a mechanical system on the group of rigid body motions $SE(3)$. The expressions explicitly demonstrate that, when a mechanical system is represented intrinsically using the Levi-Civita connection, the system can be naturally interpreted as a double integrator.

\subsubsection{Example-1: Unconstrained Mechanical Systems on the Circle}\label{Secn:MechSysOnS}
In this section we consider the special class of mechanical systems that evolve on the Lie-group $\mathbb{S}$.
The space of all possible vector fields on $\mathbb{S}$, referred to as the tangent bundle to $\mathbb{S}$ is denoted by $T\,\mathbb{S}\equiv \mathbb{S}\times \mathbb{R}$. The kinetic energy of the system defines a Riemannian metric on $\mathbb{S}$. It  is defined by $\langle\langle \zeta,\eta\rangle\rangle=\mathbb{I}(\theta)\zeta\eta$ where the inertia $\mathbb{I}(\theta)>0$ and $\zeta,\eta\in \mathbb{R}$. The unique Levi-Civita connection that corresponds to the Riemannian metric
$\langle\langle \zeta,\eta\rangle\rangle$ on $\mathbb{S}$ is explicitly given by
\begin{align*}
\nabla_\zeta \eta &=d\eta(\zeta)+\Gamma^1_{11}(\theta)\, \zeta \eta,
\end{align*}
where
\begin{align*}
\Gamma^1_{11}(\theta)&=\frac{1}{2\mathbb{I}}\dfrac{\partial \,\mathbb{I}}{\partial \theta}.
\end{align*}
The significance of the Levi-Civita connection is that it satisfies the metricity condition given by
\begin{align*}
\mathcal{L}_\xi\langle\langle \zeta,\eta\rangle\rangle&=\langle\mathbb{I}\nabla_\xi \zeta,\eta\rangle+\langle\mathbb{I}\nabla_\xi \eta,\zeta\rangle,
\end{align*}
for vector fields $\xi(\theta),\zeta(\theta),\eta(\theta)\in T\mathbb{S}$.
A mechanical system on $\mathbb{S}$ with kinetic energy equal to $\frac{1}{2}\langle\langle \omega,\omega\rangle\rangle=\frac{1}{2}\mathbb{I}(\theta)\omega^2$ is then intrinsically represented by
\begin{align}
\dot{\theta}&=\omega,\label{eq:KinematicsS}\\
\mathbb{I}\,\nabla_\omega\omega &=\tau\label{eq:MechanicalSystemS}
\end{align}
where $\tau\in \mathcal{R}$ is the generalized force. Explicitly we have that (\ref{eq:MechanicalSystemS}) is given by 
$\mathbb{I}\,\dot{\omega}+ \mathbb{I}\,\Gamma^1_{11}(\theta)\, \omega^2=\tau$.

\subsubsection{Example-2: Unconstrained Rigid Body Motion}\label{Secn:RigidBodyMotion}
The configuration space of rigid body motion is the Euclidian group $SE(3)=\mathbb{R}^3\times SO(3)$. 
Denote by $g(t)\triangleq(o(t),R(t))\in SE(3)$ the configuration of the rigid body where $o(t)\in \mathbb{R}^3$ is the position of the center of mass of the rigid body with respect to some inertial frame $\mathbf{e}$ and $R(t)\in SO(3)$ is the orientation of the rigid body that relates a body fixed frame $\mathbf{b}$ (fixed at the center of mass) to the inertial frame $\mathbf{e}$ by the relationship $\mathbf{b}=\mathbf{e}R(t)$. 

The body angular velocity of the body, $\Omega\in \mathbb{R}^3$, is given by $R^T\dot{R}=\widehat{\Omega}\in so(3)$ where $\widehat{\Omega}\in so(3)$ is the skew symmetric version of $\Omega$ that is given by the isomorphism $\mbox{}\widehat{\mbox{}}: \mathbb{R}^3\mapsto so(3)$, explicitly expressed by
\begin{align*}
\Omega=\begin{bmatrix}
\Omega_1\\
\Omega_2\\
\Omega_3 
\end{bmatrix}
\mapsto
\widehat{\Omega}=\begin{bmatrix}
0 & -\Omega_3 & \Omega_2\\
\Omega_3 & 0 & -\Omega_1\\
-\Omega_2 & \Omega_1 &0
\end{bmatrix}.
\end{align*}
We will use the following notations to represent the canonical basis elements of $\mathbb{R}^3$:  $e_1=\begin{bmatrix} 1 & 0 &0 \end{bmatrix}^T$, $e_2=\begin{bmatrix} 0 &1 & 0 \end{bmatrix}^T$, and $e_2=\begin{bmatrix} 0 & 0 & 1\end{bmatrix}^T$.

The spatial angular velocity of the body is given by
$\omega=R\Omega$. We will use the convention of using upper case greek letters to denote body variables while using the corresponding lower case letters to denote its spatial version. The two different trivializations of $TSO(3)$, arising due to left and right translations, give rise to two different representations of $\dot{g}$ explicitly given by either $(V,{\Omega})$ where $V=R^T\dot{o}$ or $(\dot{o},{\omega})$ respectively.

Consider a constant symmetric positive definite (semi-definite) matrix $\mathbb{I}_\nu$. This induces a left-invariant Riemannian (singular Riemannian) metric on
$SO(3)$ by the relationship $\langle\langle R\widehat{\Omega},R\widehat{\Psi}\rangle\rangle\triangleq\langle \mathbb{I}_\nu\Omega,\Psi\rangle=\Psi\cdot \mathbb{I}_\nu\Omega$. If $\omega\triangleq R\Omega$, $\psi\triangleq R\Psi$ then expressing the left-invariant metric using spatial angular velocities we have $\langle\langle\widehat{\omega}R,\widehat{\psi}R\rangle\rangle= \langle\mathbb{I}_\nu^R\omega,\psi\rangle=\psi\cdot\mathbb{I}_\nu^R\omega$ where now $\mathbb{I}_\nu^R\triangleq R\mathbb{I}_\nu R^T$ represents the spatial version of the left-invariant inertia tensor. On the other hand it induces a right-invariant Riemannian (singular Riemannian) metric on
$SO(3)$ by the relationship $\langle\langle \widehat{\omega}R,\widehat{\psi}R\rangle\rangle\triangleq\langle \mathbb{I}_\nu\omega,\psi\rangle=\langle \mathbb{I}_\nu^L\Omega,\Psi\rangle$ where now $\mathbb{I}_\nu^L=R^T\mathbb{I}_\nu R$ represents the body version of the right-invariant inertia tensor. A metric that is both left- and right-invariant is called a bi-invariant metric. It is clear that if $\mathbb{I}_\nu$ is bi-invariant then it must be of the form $\mu\,I_{3\times 3}$ for some positive $\mu$.
At times, for brevity, we will omit the superscript $R$ or $L$ where the version being used will be clear from the angular velocity variables being used.

The two different velocity representations of rigid body rotations arising due to the use of left or right translations on the group will give two different expretions for the Levi-Civita connection corresponding to a given left-invariant Riemannian metric on $SO(3)$.  That is, the Levi-Civita connection will have two expressions depending on whether one uses body angular velocities or spatial angular velocities to compute them. For further details about these computations we refer to the paper \cite{MaithripalaAutomatica}. Specifically it can be shown, using the Koszul formula, that for left-invariant metrics induced by $\mathbb{I}_\nu$,
\[
\mathbb{I}_\nu\nabla^\nu_\xi \eta \triangleq \mathbb{I}_\nu d{\eta}(\xi)+\frac{1}{2}\left(\pm\,\mathbb{I}_\nu(\xi\times  \eta)- \left({\mathbb{I}_\nu\eta}\times {\xi} +{\mathbb{I}_\nu\xi}\times{\eta}  \right)\right),
\]
and for right-invariant metrics induced by $\mathbb{I}_\nu$,
\[
\mathbb{I}_\nu\nabla^\nu_\xi \eta \triangleq \mathbb{I}_\nu d{\eta}(\xi)+\frac{1}{2}\left(\pm\,\mathbb{I}_\nu(\xi\times  \eta)+\left({\mathbb{I}_\nu\eta}\times {\xi} +{\mathbb{I}_\nu\xi}\times{\eta}  \right)\right),
\]
and for any bi-invariant metric
\[
\nabla^{bi}_\xi \eta \triangleq d{\eta}(\xi)\pm\frac{1}{2}\,\xi\times  \eta.
\]
In these expressions when $\xi,\eta$ are body angular velocities the $+$ needs to be chosen while when $\xi,\eta$ are spatial angular velocities the $-$ needs to be chosen. Note that these expressions make sense even when $\mathbb{I}_\nu$ is positive semi-definite and hence for singular Riemannian metrics as well.

Armed with these notations the unique Levi-Civita connection corresponding to the Riemannian metric
$\langle\langle (v,\zeta),(u,\eta)\rangle\rangle\triangleq \left(mv\cdot u+\zeta \cdot \mathbb{I}_b\eta\right)$ on $SE(3)$ that results from the kinetic energy of rigid body motions is given by
\begin{align*}
\mathbb{I}\nabla^{\tiny se}_{(v,\zeta)}{(u,\eta)}&=
\begin{bmatrix}
m\,d{u}(v)\\ \mathbb{I}_b\nabla^b_{\zeta}\eta
\end{bmatrix}.
\end{align*}
Then the intrinsic acceleration denoted by $\nabla^{\tiny se}_{\dot{g}}{\dot{g}}$ on $SE(3)$ is defined by
\begin{align*}
\mathbb{I}\nabla^{\tiny se}_{\dot{g}}{\dot{g}}&=
\begin{bmatrix}
m\,\ddot{o}\\ {\mathbb{I}_b}\dot{\Omega}-\left(\mathbb{I}_b\Omega\times \Omega\right)
\end{bmatrix}
=
\begin{bmatrix}
m\,\ddot{o}\\ {\mathbb{I}_b^R}\dot{\omega}-\left(\mathbb{I}_b^R\omega\times \omega\right)
\end{bmatrix},
\end{align*}
and thus a mechanical system on $SE(3)$ takes the intrinsic form
\begin{align*}
\mathbb{I}\nabla^{\tiny se}_{\dot{g}}{\dot{g}}&=\begin{bmatrix}f \\ \tau\end{bmatrix}
\end{align*}
where $f$ is the resultant force acting on the rigid body and $\tau$ is the force moments acting on the body about the center of mass expressed in the body frame.


\section{PID control on Lie Groups}\label{Secn:PIDonLieGroups}
For completeness we will begin this section with a review of the basic underlying notions of intrinsic PD control design for configuration trajectory tracking for fully actuated mechanical systems that was developed in \cite{Koditschek,BulloTracking,MaithripalaLieGroup,MaithripalaSIAM}. The extension to include integral control for configuration trajectory tracking of fully actuated mechanical systems will be presented in Section-\ref{Secn:PIDfullyActuated} and then it will be extended to a class of inter connected under actuated mechanical systems in Section-\ref{Secn:NLPID} and finally to constrained mechanical systems on a Lie group in Section-\ref{Secn:PID4Constrained} .

\subsection{PID Control for Fully Actuated Systems}\label{Secn:PIDfullyActuated}

The advantage that one has when working with Lie groups is two fold: a.) first the group structure allows the definition of a globally defined configuration tracking error, while b.) the trivialization of $TG$ allows the definition of a globally defined velocity tracking error without having to resort to the computationally unwieldy parallel transport maps. Left and right translations provide two distinct but equivalent ways of doing this.
For a given twice differentiable configuration reference $g_r(t)$ let the corresponding left-velocity reference be given by $\zeta^L_r= g_r^{-1}\cdot\dot{g}_r $ while $\zeta^L=g^{-1}\cdot\dot{g}$ be the left-velocity of the system.
Define the left-invariant configuration error $E^L=g_r^{-1} g$ 
and the associated intrinsic velocity error $\zeta^L_E \triangleq ({E^L})^{-1}\cdot \dot{E}^L= (\zeta^L-\mathrm{Ad}_{(E^L)^{-1}}\zeta^L_r)$ \cite{Koditschek,BulloTracking,MaithripalaLieGroup,MaithripalaSIAM}.
Similarly when the velocities of interest are the right-velocities, $\zeta^R_r=\dot{g}_r\cdot g_r^{-1}$ and $\zeta^R=\dot{g}\cdot g^{-1}$, define the right error system
by using the right-invariant error $E^R=g g_r^{-1}$ and its associated right-velocity error $\zeta^R_E \triangleq \dot{E}^R\cdot (E^R)^{-1}=(\zeta^R-\mathrm{Ad}_{E^R}\zeta^R_r)$. 
For both these types of errors the closed-loop  second order error dynamics evolving on $G\times \mathcal{G}$ take the form
\begin{align}
\mathbb{I}\nabla_{\zeta_E} \zeta_E &= f_u+\Delta_d-f_r(E,\zeta_E,\eta_r),  \label{eq:ErrorLieGorup0}
\end{align}
with
\begin{align}
f_r&=\mathbb{I}\nabla_{\zeta_E}{\eta_r}+\mathbb{I}\nabla_{\eta_r}\zeta_E+\mathbb{I}\nabla_{\eta_r}\eta_r,\label{eq:FFterm}
\end{align}
were for left-velocity representations 
$\eta_r= \mathrm{Ad}_{(E^L)^{-1}}\zeta^L_r$ and $\nabla_{\zeta_E} \zeta_E=\nabla^L_{\zeta^L_E} \zeta^L_E$ while for right-velocity representations $\eta_r= \mathrm{Ad}_{E^R}\zeta^R_r$ and $\nabla_{\zeta_E} \zeta_E=\nabla^R_{\zeta^R_E} \zeta^R_E$.

For clarity we will first investigate the convergence properties of the trajectories of the closed loop second order error system (\ref{eq:ErrorLieGorup0}) for intrinsic PD control when uncertainties are present. 
Let $V : G \rightarrow \mathbb{R}$ be a polar Morse function on $G$ with a unique minimum at $e$. Without loss of generality we assume that $V(e)=0$. It is known that polar Morse functions are good candidate local distance functions in a neighborhood of its global minimum \cite{Koditschek}.
Thus in this series of work we use the error function to quantify the configuration tracking error. We will say that the tracking error is small if the corresponding error function value is small. 

Let $\mathrm{grad}\,V(E)$ be the gradient of $V(E)$ at $E$. The gradient at $E$ is defined by the relationship
$
\left.\frac{d}{ds}\right |_{s=0}V(E(s))=\langle\,dV, v_E\rangle =\langle\langle\mathrm{grad}\,V(E(s)), v_{E}\rangle \rangle_G,
$ where $E(s)$ is a locally defined curve such that $E(0)=E$ and $
\left.\frac{d}{ds}\right |_{s=0}E=v_E$
for all $v_E\in T_E{G}$. Note that this definition also implies that $dV=\mathbb{I}\,\mathrm{grad}V$.
Depending on the choice of representation of the error system define $\eta_E\in \mathcal{G}$ to be either
 $\eta^L_E=(E^L)^{-1} \cdot \mathrm{grad}\,V(E^L)$ or
$\eta^R_E= (\mathrm{grad}\,V(E^R))\cdot (E^R)^{-1}$ respectively.
Consider the intrinsic nonlinear PD controller \cite{Koditschek,BulloTracking,MaithripalaLieGroup,MaithripalaSIAM}
\begin{align}
{f_u}&=-\mathbb{I}(k_p\eta_E+k_d \zeta_E)+f_r.\label{eq:PDcontrol}
\end{align}
The proportional control action is given by the first term, the second term gives the derivative action and the last term is the feedforward action.

In order to investigate the convergence properties of the trajectories of the error system define the positive definite function $W: G\times\mathcal{G} \to \mathbb{R}$ by
$W=k_pV(E)+\frac{1}{2}\langle\langle \zeta_E,\zeta_E \rangle \rangle$ where $\langle\langle \cdot,\cdot \rangle \rangle$ denotes either $\mathbb{I}^L_g$ or $\mathbb{I}^R_g$ depending on whether the error is left-invariant or right-invariant respectively.
Taking the derivative of $W$ along the trajectories of (\ref{eq:ErrorLieGorup0}) one obtains
\begin{align*}
\dot{W} &= k_p \langle\langle \eta_E,\zeta_E \rangle \rangle + \langle ({f_u}+\Delta_d-f_r),\zeta_E\rangle
= -k_d \langle\langle \zeta_E,\zeta_E \rangle \rangle + \langle\langle \Delta_d+\Delta_\epsilon,\zeta_E \rangle \rangle
\leq -(k_d||\zeta_E||-||\Delta_d+\Delta_\epsilon||)||\zeta_E||.
\end{align*}
Here $\Delta_\epsilon$ represents the effects due to parametric uncertainty in the inertia and the actuation models.
This shows that if the Lie group $G$ is compact and if the unmodelled forces and parametric uncertainties represented by $\Delta_d+\Delta_\epsilon$ are bounded then the trajectories of the closed loop remain globally bounded.
When $\Delta_d+\Delta_\epsilon$ is zero then LaSalle's invariance principle implies that the trajectories converge to the largest invariant set contained in $\dot{W}\equiv 0$. This set contains points of the form $(E,\zeta_E)=(\bar{E},0)$ where $\bar{E}$ satisfy $\eta_E(\bar{E})=0$. The points $\bar{E}$ are the critical points of the error function and $(\bar{E},0)$ are the equilibria of the system when $\Delta_d+\Delta_\epsilon=0$. Since $V(E)$ is a polar Morse function there are only finitely many critical points with one being a global minimum. Thus the largest invariant set is made-up of only a finite number of equilibrium points, $(\bar{E}_i,0)$, out of which only the  one corresponding to the global minimum of $V(E)$ is almost-globally stable \cite{Koditschek,BulloTracking,MaithripalaLieGroup}. 
Thus in the absence of disturbances and parametric uncertainty, the intrinsic PD controller (\ref{eq:PDcontrol}) ensures that $\lim_{t \to \infty}{E(t)} \to e$ for all initial conditions except those from a set of measure zero, corresponding to the local maxima and saddles of the polar Morse function $V(E)$ and the stable manifolds of the saddles.

This argument also shows that in the presence of bounded constant disturbances and modeling errors in the parameters and control moments, the trajectories of the error system can be made to converge almost globally to an arbitrarily small neighborhood of $(e,0)$ by picking $k_p,k_d$ to be sufficiently large. The size of the neighborhoods of $(e,0)$ are characterized by the level sets of $W$. However, picking $k_p, k_d$ large has the undesirable effect of magnifying measurement noise. Furthermore, even small constant attitude errors may cause unacceptable behavior such as translational drift in vehicles. As for linear systems, adding integral action guarantees that the attitude errors will go to zero, and avoids the use of large gains.

The intrinsic PID controller we have proposed in \cite{MaithripalaAutomatica} is
\begin{align}
\mathbb{I}\nabla_{\zeta_E} \zeta_I &= \mathbb{I}\eta_E,\label{eq:IntrinsicI}\\
f_u &= -\mathbb{I}(k_p\eta_E+k_d\zeta_E+k_I\zeta_I)+f_r(E,\zeta_E,\eta_r),\label{eq:PID}
\end{align}
where $\zeta_I\in  \mathcal{G}$.
The second order system (\ref{eq:ErrorLieGorup0}) together with the integrator (\ref{eq:IntrinsicI}) define a dynamic system on $\mathcal{X}\triangleq G\times \mathcal{G}\times\mathcal{G}$. 
Equation (\ref{eq:IntrinsicI}) essentially defines the integral term $\zeta_I$. That equation says that the gradient of $V(E)$ should be the covariant derivative of the integral term along the velocity of the error system. One could say loosely that the term $\zeta_I$ is the ``covariant integral'' of the gradient of $V(E)$. When the controller variables correspond to a left-velocity system we will call the controller a \textit{left-PID controller} while when the variables correspond to a right-velocity system we will call the controller a \textit{right-PID controller}.

\begin{remark}
The Euclidean space $\mathbb{R}^n$ is a Lie group, with the group operation of addition. On this Lie group, for mechanical systems with constant mass matrix (constant metric), $M$, the control  (\ref{eq:IntrinsicI})--(\ref{eq:PID}) with the choice of polar Morse function $V(E)=E^TE/2$ results in the well known expressions $\dot{E}_I=E$, $f_u=-M(k_p E +k_d \dot{E}+k_IE_I-\ddot{x}_r)$, where $x_r\in \mathbb{R}^n$ is the reference and $E=(x-x_r)$.
\end{remark}

Analogous to the linear case, in order to investigate the convergence properties of the closed loop system we consider the following closed loop error dynamics defined on $G\times \mathcal{G}\times\mathcal{G}$:
\begin{align}
\mathbb{I}\nabla_{\zeta_E} \zeta_E &= -\mathbb{I}(k_p\eta_E+k_d\zeta_E+k_I\zeta_I)+\Delta_d +\Delta_\epsilon \label{eq:ErrorLieGorup}\\
\mathbb{I}\nabla_{\zeta_E} \zeta_I &= \mathbb{I}\eta_E,\label{eq:IntrinsicI2}
\end{align}
where the parameter error term $\Delta_\epsilon$ is given by
\begin{align*}
\Delta_\epsilon &=-\epsilon(k_p\eta_E+k_d\zeta_E+k_I\zeta_I)+\epsilon_{f_r}.
\end{align*}
Here $\epsilon$ is a linear operator, $\epsilon : \mathcal{G} \to \mathcal{G}^*$ and $\epsilon_{f_r}\in \mathcal{G}^*$ that depends on the inertia uncertainties $(\mathbb{I}^{-1}{\mathbb{I}}_0-I_{n\times n})$ and the actuator model uncertainties where $\mathbb{I}$ is the actual inertia tensor and ${\mathbb{I}}_0$ is the nominal inertia tensor with $n=\mathrm{dim}(G)$.
It is useful to note that when the reference velocity $\zeta_r$ and the disturbance $\Delta_d$ are constant, the equilibria of the error dynamics are of the form
\[
(\bar{E}_i,0,(\mathbb{I}+\epsilon)^{-1}({\Delta}_d+\epsilon_{f_r})/k_I)
\] 
where $\bar{E}_i$ are the critical points of $V(E)$. Similar to the case of the PD controller the convergence properties of the error system trajectories completely depend on the assigned potential energy $V(E)$.  All these results explicitly assume:
\begin{assumption}\label{AssumptionMain}
The reference trajectories are twice differentiable and bounded.
The disturbance and unmodelled forces represented by $\Delta_d(t)$, and the effect of the parametric uncertainties of the system represented by $\Delta_\epsilon$ are bounded.  
\end{assumption} 

The function, its gradient and the Hessian of a polar Morse function $V(E)$ are bounded  on a compact set \cite{Koditschek,Morse}. 
Let $\mu$ denote the bound on the Hessian and let $\lambda$ be defined by
$\lambda =\sup_{\mathcal{X}}\frac{\langle\langle \eta_E,\eta_E\rangle\rangle}{2 V(E)}$.
The existence of a value for $\lambda$ satisfying this inequality is guaranteed by the fact that $V(E)$ is a polar Morse function. Boundedness of the gradient implies that ${\langle\langle \eta_E,\eta_E\rangle\rangle}/{(2 V(E))}$  is finite away from points where $V(E)=0$. Since polar Morse functions have only a unique minimum and we have chosen the unique minimum value to  be zero what one needs to ensure is that the limit exists at the unique minimum. This follows from the fact that a Morse function is locally quadratic and positive definite around a  
minimum \cite{Morse}. In fact what the above expression means is that the function $V(E)$ is quadratically bounded from below by $\langle\langle \eta_E,\eta_E\rangle\rangle/(2\lambda)$ at the global minimum. A fact that will ensure locally exponential convergence.

Let the controller gains $k_p,k_I,k_d>0$ be chosen to satisfy the following inequalities involving $\lambda$ and $\mu$. 
\begin{align}
0<k_I<\frac{k_d^3(1-\delta^2)}{\mu},\label{eq:kICond}\\
 k_p>\max \left\{k_1,k_2,2\kappa k_d^2\right\},\label{eq:kpCond}
 \end{align}
 where,
 \begin{align*}
 k_1&=\frac{k_I}{2k_d}\left(\sqrt{1+\frac{16\lambda\kappa^2k_d^2}{k_I}}-1\right),\\
 k_2&=\frac{ \lambda k_I^2}{2k_d^4}\left({1 +\sqrt{1+\frac{4k_d^3 (k_I^2+4\kappa k_d^3(1+\kappa k_d^3))}{\lambda k_I^3}}}\right),
\end{align*}
where $0<\kappa <2/\mu$ and $\delta \triangleq |(\kappa \mu-1)|<1$.

The following corollary that was proven in \cite{MaithripalaAutomatica} also follows from the Theorem-\ref{theom:Theorem1} that was proven in \cite{RollingHoopACC2017}. For the sake of completeness we reproduce the proof of this theorem in the appendix.
\begin{corollary}\label{CorAutomatica}
Let Assumption-\ref{AssumptionMain} hold and $\mathcal{X}\subset G\times \mathcal{G}\times\mathcal{G}$  be a compact sub set that contains $(e,0,0)$.
If $V(E)$ is a polar Morse function on $G$ with a unique minimum at $e$ and the gains  $k_I$, $k_d$, and $k_p$ of the intrinsic PID controller (\ref{eq:IntrinsicI})--(\ref{eq:PID}) are chosen such that they satisfy (\ref{eq:kICond})--(\ref{eq:kpCond}),
then the followings hold:
\begin{enumerate}
\item the tracking error quantified by $V(E)$, the velocity error $\zeta_E$, and the integrator state $\zeta_I$ are bounded for all initial conditions in $\mathcal{X}$,
\item  by picking sufficiently large gains the asymptotic tracking error quantified by $V(E)$, the velocity error, $\zeta_E$, and the integrator state, $\zeta_I$, can be made arbitrarily small for almost all initial conditions in $\mathcal{X}$,
\item If $\Delta_d$ and $\zeta_r$ are constant then $\lim_{t \to \infty}g(t) = g_r(t)$ locally exponentially  for almost all initial conditions in $\mathcal{X}$.
\item In the absence of uncertainty $\lim_{t \to \infty}g(t) = g_r(t)$ locally exponentially  for almost all initial conditions in $\mathcal{X}$.
\end{enumerate}
If the Lie group $G$ is compact then $\mathcal{X} = G\times \mathcal{G}\times\mathcal{G}$ is the entire state space.
\end{corollary}

To the best of our knowledge, Corollary-\ref{CorAutomatica} is the strongest such result reported in the literature for robust smooth state feedback tracking for fully actuated mechanical systems on a general Lie group including the widely treated Lie group of rigid body rotations 
$\mathrm{SO}(3)$. In the case of $\mathrm{SO}(3)$ one can show that the convergence is guaranteed for all initial conditions except three initial conditions that correspond to the unstable equilibria of the closed loop system. 

\begin{remark}
We highlight that the conditions 1--4 stated in Corollary-\ref{CorAutomatica} have a one-to-one correspondence with the properties of the linear PID controller for the double integrator system on $R^{n}$ where global exponential stability in the linear case must be replaced with almost-global and locally exponential stability for simple mechanical systems on Lie groups.
\end{remark}


\subsection{PID Control for a class of Underactuated Mechanical Systems}\label{Secn:NLPID} 

In the following we consider a class of interconnected mechanical systems on Riemannian manifolds. Each subsystem evolves on a configuration space $G_\nu$ and has an inertia tensor $\mathbb{I}_\nu$. Here $\nu$ is either $s$ or $a$, denoting the corresponding subsystem of the interconnected system. Denote by $\dot{g}_{\nu}=v_\nu$ and by $\nabla^\nu$ be the unique Levi-Civita connection corresponding to $\mathbb{I}_\nu$.

The class we consider is of the form
\begin{align}
\mathbb{I}_s\nabla^s_{v_s}{v_s} &=\Delta_s+\tau_{s}(v_s,v_a)+\tau_u,\label{eq:System_s}\\
\mathbb{I}_a\nabla^a_{v_a}{v_a} &=\tau_V(g_a)+\Delta_a+\tau_{a}(v_s,v_a)+B\tau_u,\label{eq:System_a}
\end{align}
where, $\Delta_s$ and $\Delta_a$ represent unmodelled forces and disturbances,  $\tau_V$ is a conservative force, and $\tau_{s}$ and $\tau_{a}$ are velocity dependent interaction forces. 
The velocity interaction forces $\tau_{\nu} \::\: T_{g_s}{G}_s\times T_{g_a}{G}_a\mapsto T_{g_\nu}{G}^*_\nu$ are assumed to be bilinear in the arguments. The system denoted by $s$ will be referred to as the \textit{ output system} and the system denoted by $a$ will be referred to as the \textit{actuation system}.

Both subsystems are assumed to be have configuration spaces of the same dimension and be fully actuated with respect to the common controls $\tau_u$. These assumptions do not impose any significant physical restrictions since for any interconnected mechanical system Newton's $3^{\mathrm{rd}}$ law imposes that all interactions are equal and opposite and that the gyroscopic forces are quadratic in the velocities.

In this section we consider the problem of ensuring output tracking of the interconnected mechanical system (\ref{eq:System_s})--(\ref{eq:System_a}).
We assume that the output of the system takes values on a Lie-Group $G_y$ and is only a function of the configuration of the output subsystem denoted by `$s$'. That is the output is a smooth function
\begin{align}
y\:\:\: : \:\:\: G_s \mapsto G_y.\label{eq:Output}
\end{align}
We also assume that the output is relative degree two with respect to the input $\tau_u$ and that the potential function $V_a(g_a)$, that gives rise to the conservative force $\tau_V(g_a)$ of the auxiliary system, is a Polar Morse function.

The control problem that we solve in this section is that of ensuring the almost-semi-global and local exponential convergence of the smooth output $y$ to $e_y$ the (the identity of $G_y$) while ensuring that $\lim_{t\to \infty} v_s(t)=0$ and $\lim_{t\to \infty} v_s(t)=0$. 
Let $V_y : G_y \mapsto \mathbb{R}$ be a polar Morse function on $G_y$ with a unique minimum at the identity, $e_y$, of $G_y$. Let $v_s(g_s)=V_y(y(g_y))$ and $\eta_e\in T_{g_s}G_s$ be the gradient of $V_s$ that is defined by $\langle\langle\eta_e,v_s\rangle\rangle=dV_s(v_s)$.

Consider the nonlinear PID controller
\begin{align}
&\mathbb{I}_s\nabla^s_{v_s} {v_I} =\mathbb{I}_s{\eta_e},\label{eq:GeneralI}\\
&\tau_u=-\mathbb{I}_s(k_p{\eta_e}+k_d{v_s}+k_I{v_I})-k_cB^{-1}\mathbb{I}_av_a,\label{eq:GeneralPID}
\end{align}
where $v_I,\eta_e\in T_{g_s}G_s$

In \cite{RollingHoopACC2017} we will prove the following theorem. For the sake of completeness we also reproduce it here the Appendix.
\begin{theorem}\label{theom:Theorem1} 
If the conditions of Assumption-\ref{AssumptionMain} hold and
if the gains of the nonlinear PID controller (\ref{eq:GeneralI})--(\ref{eq:GeneralPID}) are chosen to satisfy (\ref{eq:kICond}) and (\ref{eq:kpCond})
then for almost all $(g_s(0),{v_s}(0),{v_I}(0),g_a(0),{v_a}(0))\in \mathcal{X}$, where $\mathcal{X}$ is compact, the output $\lim_{t\to \infty} y(g_s(t))=e_y$ and $\lim_{t\to \infty} ({v_s}(t),{v_a}(t))=(0,0)$ semi-globally and locally exponentially in the presence of bounded constant disturbances and bounded parameter uncertainty. Here  $1/\mu<\kappa<2/\mu$ where $\mu=\mathrm{max}\,||\nabla^s{\eta_e}||$ on $\mathcal{X}$ restricted to $G_s$ and $\lambda=\max_{{\mathcal{X}|_{G_s}}}\frac{\langle\langle\eta_e,\eta_e\rangle\rangle_s}{2V_s(g_s)}$.
\end{theorem}

From the proof of this theorem it can be shown that the following Corollary and the Corollary-\ref{CorAutomatica} that was proven in \cite{MaithripalaAutomatica} hold. 
\begin{corollary}
If the velocity interconnection terms $\tau_a\equiv 0$ and $\tau_b\equiv 0$ and the Lie-groups $G_s$, $G_a$, and $G_y$ are compact then the convergence is almost-global and locally exponential.
\end{corollary}


\subsection{PID Control for Constrained Mechanical Systems}\label{Secn:PID4Constrained}

In here we extend the PID controller presented in Section-\ref{Secn:PIDfullyActuated} for fully actuated unconstrained mechanical systems to the setting of output tracking of constrained mechanical systems on Lie groups. Specifically we consider a mechanical system with constraints described by a smooth nonsingular $m$-dimensional distribution $\mathcal{D}$, control forces $\gamma$. 
We have shown in Section-(\ref{Secn:CovDerivative}) that such a system is described by the expression (\ref{eq:ConstrainedMech}). We assume that the output of the system
$y : G \mapsto \mathcal{Y}$ is a smooth onto map from the configuration space to a smooth manifold $\mathcal{Y}$ and that the dimension of $\mathrm{Span}(P_{\mathcal{D}^*}(\gamma))$ is $m$. The problem that we consider is ensuring that $\lim_{t\to \infty} y(t)=y_{\mathrm{ref}}(t)$ for some twice differentiable reference output $y_{\mathrm{ref}}(t)$. Let $V_y : \mathcal{Y}\times \mathcal{Y}\to \mathbb{R}$ be a distance function on $\mathcal{Y}$. We will assume that the error function $V_{y_0} :\mathcal{Y}\mapsto \mathbb{R}$ explicitly given by $V_{y_0}(\cdot)=V_y(\cdot,y_0)$ is a polar morse function for every fixed $y_0\in \mathcal{Y}$. 
We propose the following PID controller for solving this tracking problem:
\begin{align}
\mathbb{I}\nabla_{\dot{g}}v_g^I&=-(\nabla_{\dot{g}}P_{\mathcal{D}^*_c})(\mathbb{I}v_g^I)+P_{\mathcal{D}^*}(dV_{y_0}),\label{eq:ConstrianedI}\\
P_{\mathcal{D}^*}(\gamma)&=-P_{\mathcal{D}^*}(dV_{y_0})-k_dP_{\mathcal{D}^*}(\mathbb{I}\dot{g})-k_IP_{\mathcal{D}^*}(\mathbb{I}v_g^I).\label{eq:ConstrianedPID}
\end{align}
Here $dV_{y_0}$ is the differential of $V_{y_0}$ extended to $G$ that is of $V_{y_0}(y(g))$ and is explicitly given by $\langle dV_{y_0},v_g\rangle=\left.\dfrac{d}{ds}\right|_{s=0}V_{y_0}(y(g(s)))$ where $g(s)$ is smooth locally defined curve such that $g(0)=g$ and $\left.\dfrac{dg}{ds}\right|_{s=0}=v_g$ for all $v_g\in T_gG$.

In \cite{UdariCDC2017} we prove the following theorem.
\begin{theorem}\label{Thoem:ConstrainedPID}
If the conditions of Assumption-\ref{AssumptionMain} hold and if the gains of the nonlinear PID controller (\ref{eq:ConstrianedI})--(\ref{eq:ConstrianedPID}) are chosen to satisfy (\ref{eq:kICond}) and (\ref{eq:kpCond})
then for almost all $(g(0),\dot{g}(0),{v^I_g}(0))\in \mathcal{X}\subset $, where $\mathcal{X}$ is compact, the output $\lim_{t\to \infty} y(g(t))=y_0$ semi-globally and locally exponentially in the presence of bounded constant disturbances and bounded parameter uncertainty. Here  $1/\mu<\kappa<2/\mu$ where $\mu=\mathrm{max}\,|| \mathbb{I}^{-1}dV_{y_0}||$ on $\mathcal{X}$ restricted to $G$ and $\lambda=\max_{{\mathcal{X}|_{G}}}\frac{\langle\langle \mathbb{I}^{-1}dV_{y_0}, \mathbb{I}^{-1}dV_{y_0}\rangle\rangle}{2V_{y_0}(y(g))}$. The convergence is almost global if ${G}$ is compact.
\end{theorem}


\section{Examples of PID Control on Lie Groups}\label{Secn:ExamplesPID_LieGroups}
In this section we will provide examples of the application of the general controllers developed in
 Sections \ref{Secn:PIDfullyActuated}, \ref{Secn:NLPID} and \ref{Secn:PID4Constrained}. The PID controller developed for fully actuated mechanical systems on Lie groups will be demonstrated for the almost global locally exponential tracking of a multi-rotor-autoguided-aerial-vehicle (MRAV) in Section-\ref{Secn:AttitudeStabilization}. This section will present simulation results as well as experimental results. The extension of the PID controller to a class of under actuated interconnected mechanical systems will be demonstrated on three different physical systems: an inverted pendulum on a cart in Section-\ref{Secn:IPC}, a hoop rolling without slip on an inclined plane of unknown angle of inclination in Section-\ref{Secn:RollingHoop}, and a sphere rolling without slip on an inclined plane of unknown angle of inclination in Section-\ref{Secn:RollingSphere}. Finally the PID controller for constrained mechanical systems on Lie groups developed in \ref{Secn:PID4Constrained} is demonstrated for almost global locally exponential tracking of a spherical pendulum.

\subsection{Attitude Stabilization of a Quadrotor}\label{Secn:AttitudeStabilization}

Consider the rigid body model of the quadrotor given by
\begin{align}
\dot{R} &=  R\widehat{\Omega},\label{eq:dotRquad}\\
\mathbb{I}\dot{\Omega} &= \mathbb{I}{\Omega} \times \Omega+\mathbb{I}\Delta_T+ {{\tau}^u},\label{eq:dotOmegaquad}
\end{align} 
where $R\in SO(3)$ is the attitude or the rotation matrix of the quadrotor with respect to an inertial frame $\mathbf{e}=\{\mathbf{e}_1,\mathbf{e}_2,\mathbf{e}_3\}$ where $\mathbf{e}_2$ points towards North and $\mathbf{e}_3$ points vertically upwards, $\Omega$ is the angular velocity of the quadrotor, 
$\mathbb{I}\Delta_T$ is the generalized disturbance, $\mathbb{I}$ is the inertia tensor, and $\tau^u=\mathbb{I}T^u$ represents  the control moments generated by the thrust of the rotors. The $\Delta_T$ term includes both uncertainties in the inertia tensor and the actuator models, and the external disturbance torques and uncertainties in the center of mass location. 
The rotation matrix satisfies $R^TR=I_{3\times 3}$ and the angular velocity $\Omega$ is defined by $\widehat{\Omega}=R^T\dot{R}$ where $\widehat{\Omega}$ is the skew symmetric version of $\Omega\in \mathbb{R}$. The isomorphism between $\mathbb{R}^3$  and $so(3)$, the space of $3\times 3$ skew symmetric matrices that is given by $\:\:\widehat{\mbox{}}: \mathbb{R}^3 \to so(3)$.
Neglecting the effect of rotor flapping, the nominal body moments generated by the rotors can be expressed as
\begin{align}
\tau^u=\left[\begin{array}{cccc} 0 & lc_l & -lc_l & 0\\-lc_l & 0 & lc_l &0\\-c_d & c_d & -c_d & c_d
\end{array}\right]
\left[ \begin{array}{c}
\omega_1^2 \\ \omega_2^2 \\ \omega_3^2 \\ \omega_4^2 
\end{array}\right],\label{eq:QuadMoments}
\end{align}
and the total thrust force can be expressed as
\begin{align}
f^u&=c_l(\omega_1^2 + \omega_2^2 + \omega_3^2  + \omega_4^2),\label{eq:QuadThrust}
\end{align} 
where $l$ is the distance from the center of mass of the UAV to the center of a rotor, $c_l$ is the coefficient of lift acting on a rotor, $c_d$ is the co-efficient of drag acting on a rotor, $M$ is the total mass of the quadrotor, and $\omega_i$ is the angular velocity of the $i^{\mathrm{th}}$ rotor. 

It is shown in \cite{MaithripalaAutomatica} that the the intrinsic PID controller (\ref{eq:IntrinsicI})--(\ref{eq:PID}) that is explicitly expressed below is capable of locally exponentially and almost globally stabilizing a sufficiently smooth reference configuration trajectory $R_r(t)$ in the presence of disturbances and parametric uncertainty.
\begin{align} 
E&=R_r^TR \label{eq:TrackingError}\\
\widehat{\eta}_E&=RR_r^T-R_rR^T \label{eq:netaE}\\
\widehat{\Omega}_e&=E^T\dot{E}\label{eq:OmegaE}\\
\mathbb{I}\dot{\Omega}_I&= -\frac{1}{2}\left(\mathbb{I}(\Omega \times \Omega_I) - (\mathbb{I}\Omega_I\times \Omega+\mathbb{I}\Omega\times \Omega_I)\right)+\mathbb{I}\eta_E\label{eq:IntegratorQuad}\\
{{\tau}^u}&=-\mathbb{I}(k_p\eta_E+k_d \Omega_e+k_I \Omega_I)+f_{\mathrm{ref}}.\label{eq:PIDcontrol}
\end{align}

Here $k_p, k_d$ and $k_I$ are positive gains to be chosen through simulation and then fine-tuned experimentally. Here (\ref{eq:TrackingError}) gives the tracking error, (\ref{eq:netaE}) gives the gradient of the error function $f(E)=\mathrm{trace}(I_{3\times 3}-E)$, (\ref{eq:IntegratorQuad})
is the intrinsic integrator, and (\ref{eq:PIDcontrol}) is the nonlinear intrinsic PID control action. 
The total thrust force $f^u$ is manually set to a desired value through the throttle command transmitted through the radio controller and can be used to manually control the altitude of the quadrotor.

From (\ref{eq:QuadMoments}) and (\ref{eq:QuadThrust}) we have that the required motor speeds can be calculated by
\begin{align}
\left[ \begin{array}{c}
\omega_1^2 \\ \omega_2^2 \\ \omega_3^2 \\ \omega_4^2 
\end{array}\right]=
A^{-1}\begin{bmatrix}
f^u \\ \tau^u
\end{bmatrix}
\end{align}
where
\begin{align*}
A=\left[\begin{array}{cccc} c_l & c_l & c_l & c_l\\ 0 & lc_l & -lc_l & 0\\-lc_l & 0 & lc_l &0\\-c_d & c_d & -c_d & c_d
\end{array}\right].
\end{align*}
From this it is clear that the motor constants can also be absorbed into the unknown controller gains and thus the knowledge of these parameters are not necessary for the implementation of the controller. Therefore we stress that the implementation of the controller only requires the knowledge of the inertia tensor of the device. However, a rough estimate of these system parameters can be very helpful in finding a suitable set of gains through the use of a simulation.

\subsubsection{Simulation Verification}
We let the nominal inertia parameters be $M=0.65\,$kg, $\mathbb{I}=\mathrm{diag}\{0.004,0.004,0.006\}$ kg\,$\mathrm{m}^2$ while the actual parameters be $M_0=0.7\,$kg, $\mathbb{I}_0=\mathrm{diag}\{0.0035,0.0045,0.007\}$ kg\,$\mathrm{m}^2$. The actual actuator parameters $c_l,c_d$ were taken to be $\pm$10\% different from the nominal values. We have also assumed that the motors are not identical and that they saturate at an upper limit of 15,000 r.p.m. and at a lower limit of 2000 r.p.m. 
 
The motor speeds $\omega_i$ are determined using (\ref{eq:QuadMoments})--(\ref{eq:QuadThrust}) where we let 
$f^u=1.3g$ and $\tau^u$ given by (\ref{eq:IntegratorQuad})--(\ref{eq:PIDcontrol}). Due to the saturation of the motors we implement the controller (\ref{eq:IntegratorQuad})--(\ref{eq:PIDcontrol}) with a simple anti-windup scheme to prevent integrator windup where the the integrator update is stopped when any one of the motors saturate. We employ the almost globally and locally exponentially convergent observer proposed in \cite{MahonyHamel} to obtain the estimate $(R,\Omega)$. The observer equations are not presented here for brevity. We refers the reader to \cite{MahonyHamel,QuadICIIS2015} for these equations.

Simulations were carried out at a step size of $h=1$ ms using a Runge-Kutta numerical differentiation algorithm based on the MATLAB ODE45 function while the controller was updated every $20$ ms.
In all simulations the initial conditions of the quadrotor were chosen to correspond to an upside down stationary configuration ($R(0)=\mathrm{diag}\{1,-1,-1\}$ and $\Omega(0)=[0,0,0]\,rad/s$).  

\begin{figure}
  \centering
\begin{tabular}{cc}
\includegraphics[width=2.5in]{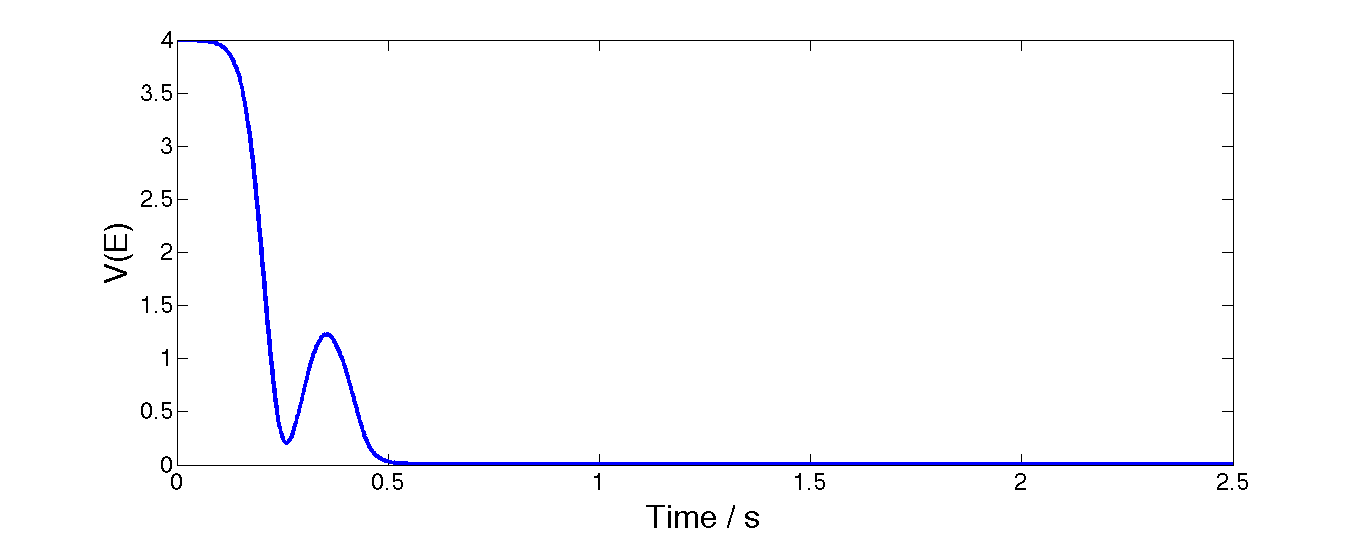} &
\includegraphics[width=2.5in]{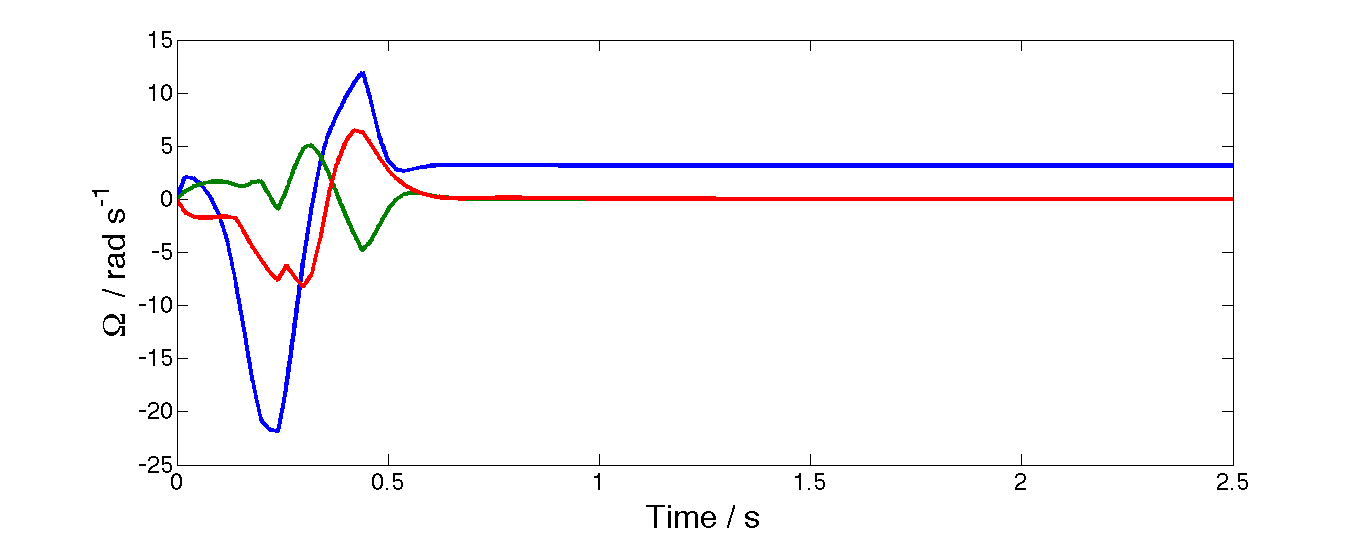}\\
\includegraphics[width=2.5in]{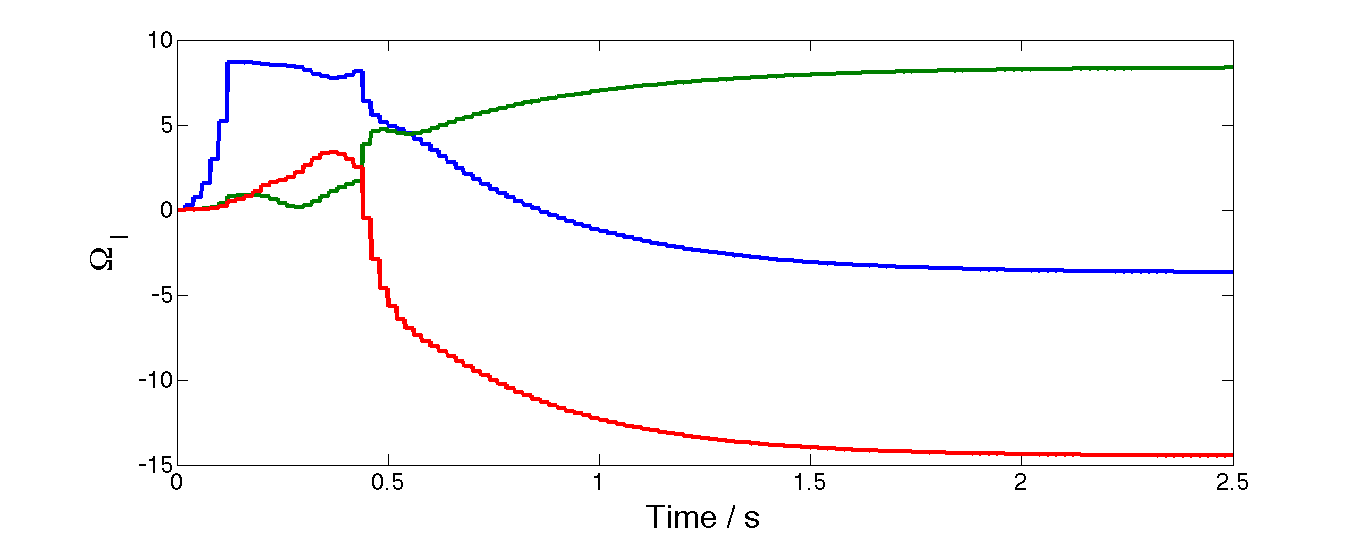} & 
\includegraphics[width=2.5in]{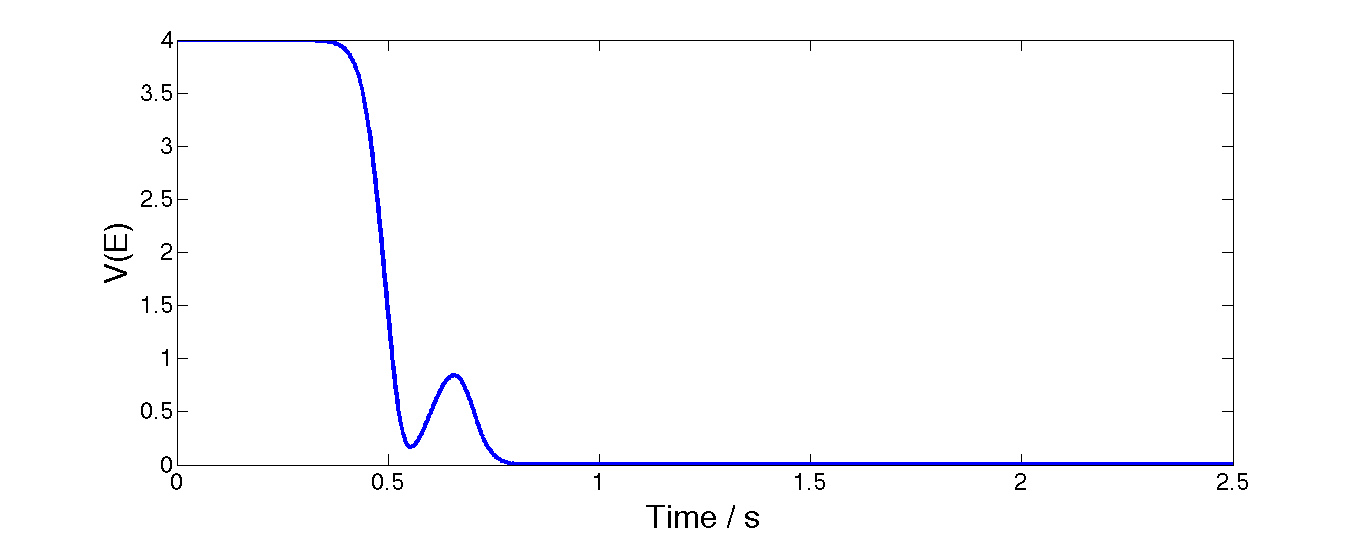}
\end{tabular}
 \caption{The tracking error $V(E)=\mathrm{trace}(I_{3\times 3}-{E})$, body angular velocities $\Omega$, the integrator state $\Omega_I$, and he motor speeds for the quadrotor tracking the attitude reference $R_r(t)=\exp{\pi te_1}$  using the control (\ref{eq:IntegratorQuad})--(\ref{eq:PIDcontrol}) with $\widehat{\tilde{\eta}}_E=\frac{1}{\mathrm{det}(\mathbb{I})}\mathbb{I}(\tilde{E}-\tilde{E}^T)\mathbb{I}$ in the presence of constant disturbances, parametric uncertainty, and actuator saturation.}
\label{Fig:TimeInvariantNoNoiseError}
\end{figure}

For the case of attitude tracking we use the control (\ref{eq:IntegratorQuad})--(\ref{eq:PIDcontrol}) along with $\hat{\tilde{\eta}}_E(t)=\frac{1}{\mathrm{det}(\mathbb{I})}\mathbb{I}(\tilde{E}-\tilde{E}^T)\mathbb{I}$and the gains $k_p=2, k_I=5, k_d=35$. Figure-\ref{Fig:TimeInvariantNoNoiseError} 
demonstrates the effectiveness of the controller when the reference is $R_r(t)=\exp{(\pi t\,e_1)}$.
In these simulations we assume that an un-modeled constant moment $\Delta_d=-g(\bar{X}\times e_3)$ is acting on the quadrotor.

\subsubsection{Experimental Verification}
The nonlinear PID controller gains were chosen to be $k_p=[100\;100\;50]$, $k_d=[30\;30\;20]$, and $k_I=[10\;10\;0]$. Low gains were used in the yaw direction to facilitate smooth take off. We stress that the controller (\ref{eq:TrackingError})--(\ref{eq:PIDcontrol}) together with the observer\cite{MahonyHamel} are entirely implemented in the onboard ATmega2560 processor and that no offboard sensing or processing is used to stabilize the quadrotor. The only external commands that are transmitted to the onboard processor are the set point command signals for yaw, pitch, roll, and throttle. These are transmitted to the quadrotor through the radio controller unit. The onboard sensor and the observer data are transmitted back to a offboard computer through the radio module for post processing purposes only. The loop time for the controller is around 5ms. It was also observed that the loop time can be delayed  to 20ms without significantly affecting the stability performance.

\begin{figure}
  \centering
\includegraphics[width=3.6in]{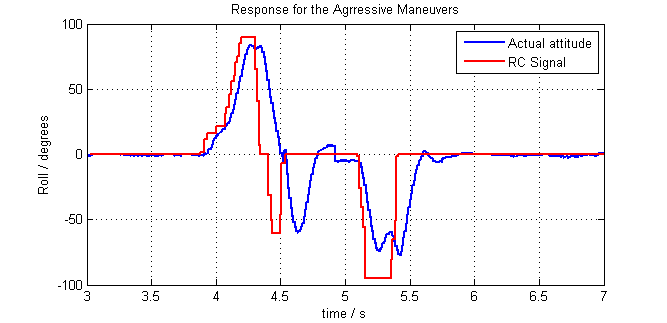}
  \caption{Variation of the roll angle of the quad rotor for a step set point command that causes a roll deviation as large as $90^0$: Blue - actual angle, Red - command given by the RC signal.}\label{Fig:stabi_Error}
\end{figure}

\begin{figure}
  \centering
\includegraphics[width=3.6in]{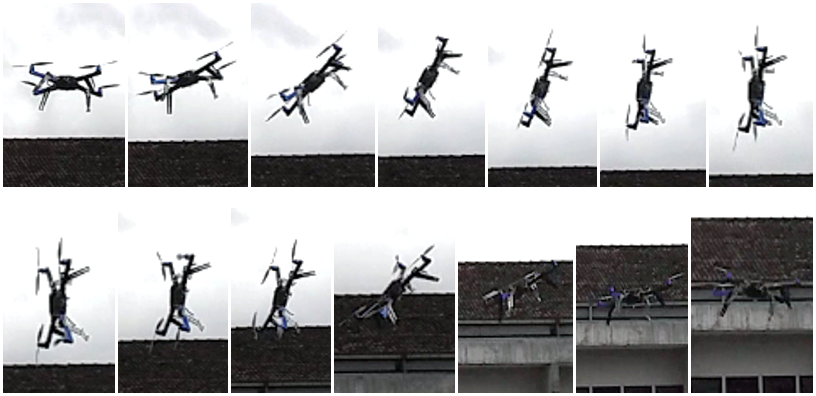}
  \caption{The actual video frame capture sequence of the quadrotor that corresponds to the results shown in Figure-\ref{Fig:stabi_Error}.}\label{Fig:capture}
\end{figure}

Figure-\ref{Fig:stabi_Error} shows excellent hovering mode stability and excellent recovery to the hovering mode for even a roll deviation as large as $90^0$.  The red curve in the figure corresponds to the roll set point command given by the radio controller and the blue curve shows the roll response of the quadrotor. The actual video frame capture sequence of the quadrotor that corresponds to the results shown in Figure-\ref{Fig:stabi_Error} is shown in figure-\ref{Fig:capture}. 

\begin{figure}
  \centering
  \begin{tabular}{cc}
\includegraphics[width=2.5in]{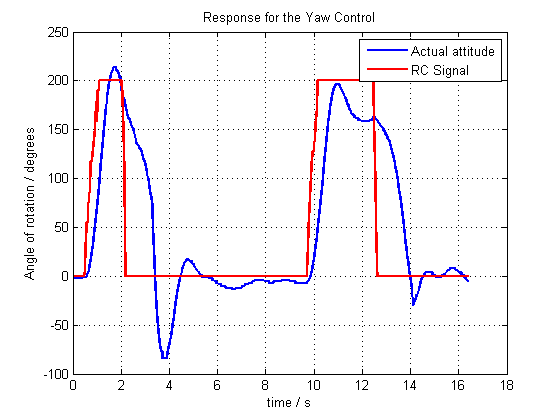}& \includegraphics[width=2.85in]{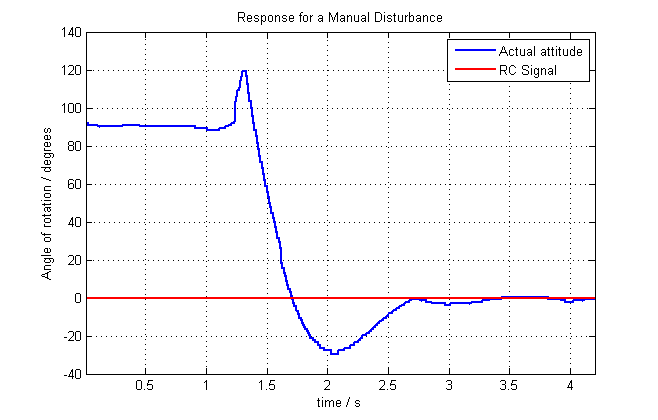}\\
(a) the yaw angle & (b) the roll angle 
\end{tabular}
  \caption{Variation of the yaw angle of the quad rotor for a yaw step set point command and the variation of the roll angle for large initial condition error.}\label{Fig:y_Error}
\end{figure}

Similarly figure-\ref{Fig:y_Error}(a) shows the test results when a yaw step set point command was transmitted through the radio controller. As seen in the figure the quadrotor is capable of recovering the hovering mode for deviations as large as $180^0$.

Finally in figure-\ref{Fig:y_Error}(b) we show the results of recovery to a specified upright configuration from large initial condition deviations.  
The quadrotor was armed and held manually at a roll angle of $90^0$ and then released 1.3 seconds later while giving an extra tilt corresponding to about $120^0$. 

\subsection{Inverted Pendulum on a Cart}\label{Secn:IPC}
\begin{figure}[ht]
	\begin{center}
		\includegraphics[width=2.5in]{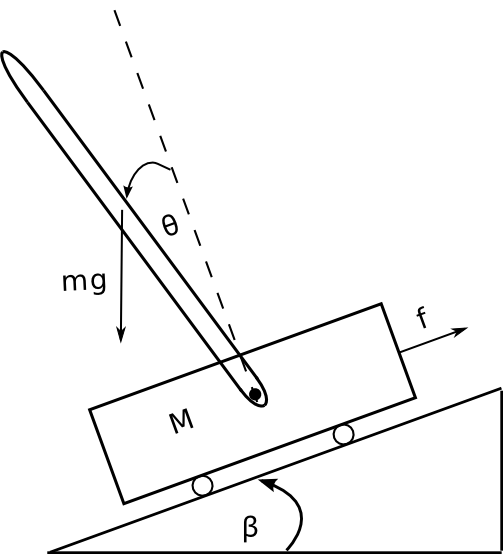}\\
		\renewcommand{\baselinestretch}{1}\selectfont
		\caption{Inverted Pendulum on a Cart.}
		\label{Fig:InvertedPendulum}
		\renewcommand{\baselinestretch}{1.5}\selectfont
	\end{center}
\end{figure}

Consider the inverted pendulum on a cart shown in Figure-\ref{Fig:InvertedPendulum}. This approximates the behavior of a segway type balancing vehicle if the yaw dynamics are relatively slow compared to the pitching dynamics. The configuration space of the system is $\mathbb{R}\times \mathbb{S}$. If $(x,\theta)\in\mathbb{R}\times \mathbb{S}$ denotes the coordinates of that represent a certain configuration of the system. Kinetic energy and potential energy of the system are given by,
\begin{align*}
KE&=\frac{1}{2}\left[M\dot{x}^2+m\left(\dot{x}^2-2L\cos\theta\dot{\theta}\dot{x}\right)+\mathbb{I}_p\dot{\theta}^2\right]\\
PE&=mgL\cos(\theta+\beta),
\end{align*}
where
$M$ is the mass of the cart, $m$ is the mass of the pendulum, $L$ is the distance to the center of mass of the pendulum from the pivot point, $\mathbb{I}_p$ is the moment of inertia of the pendulum with respect to the pivot point of the pendulum. The Euler-Lagrange equations of the system are then given by
\begin{align*}
(M+m)\ddot{x}-mL\cos\theta\ddot{\theta}+mL\sin\theta\dot{\theta}^2&=f,\\
\mathbb{I}_p\ddot{\theta}-mL\cos\theta\ddot{x}-mgL\sin(\theta+\beta)&=0.
\end{align*}
This is an underactuated mechanical system. The control problem of interest is to ensure that the output $y=(\theta(t)+\beta)\in\mathbb{S}$ is driven to zero. From the above two equations we find that 
the system evolves according to
\begin{align}
&\mathbb{I}(\theta)\ddot{\theta}+\frac{m^2L^2}{M+m}{{\dot{\theta}}}^{2}\,\sin{\theta}\cos{\theta}=
mgL\sin(\theta+\beta)+\frac{mL\cos\theta}{M+m}\,f,\label{eq:IPC}\\
&(M+m)\ddot{x}=-{\frac{mL\mathbb{I}_p}{\mathbb{I}(\theta)}{{ \dot{\theta}}}^{2}\,\sin{\theta}+\frac{{m^2}{L}^{2}g}{\mathbb{I}(\theta)}\sin{(\theta+\beta)\cos\theta}+\frac{\mathbb{I}_p}{\mathbb{I}(\theta)}\, f} ,\label{eq:xddot}
\end{align}
\noindent where $\mathbb{I}(\theta)=\left(\mathbb{I}_p-\frac{m^2L^{2}}{M+m} \cos^2 {\theta}\right)$. Let $\nabla$ be the unique Levi-Civita connection corresponding to the metric $\langle\langle\zeta,\eta\rangle\rangle=\mathbb{I}(\theta)\zeta\eta$ which is explicitly given by 
\begin{align*}
\nabla_\zeta \eta&=d\eta(\zeta)+\frac{m^2L^2\sin{2\theta}}{2\mathbb{I}(\theta)(M+m)}\zeta \eta.
\end{align*}
Then we find that (\ref{eq:IPC}) -- (\ref{eq:xddot}) can be expressed as
\begin{align*}
\dot{\theta}&=\omega,\\
\mathbb{I}(\theta)\nabla_\omega \omega&=
mgL\sin(\theta+\beta)+u.\\
v&=\dot{x},\\
(M+m)\ddot{x}&=-\frac{mL\mathbb{I}_p}{\mathbb{I}(\theta)}{{ \dot{\theta}}}^{2}\,\sin{\theta}+\frac{{m^2}{L}^{2}g}{\mathbb{I}(\theta)}\sin{(\theta+\beta)}\cos\theta+B(\theta) u
\end{align*}
where we define
\begin{align*}
u&\triangleq \frac{mL\cos\theta}{M+m}\,\,f,\:\:\:\:\:\:\:
B(\theta)\triangleq\frac{{(M+m)}\mathbb{I}_p}{(mL\cos\theta)\mathbb{I}(\theta)}
\end{align*}
The output of interest is $y=\theta \in \mathcal{S}$. Thus from the results of Section-\ref{Secn:NLPID} we have that the PID controller explicitly given by
\begin{align}
\dot{o}_I&=-\frac{m^2L^2\sin{(2\theta)}}{2\mathbb{I}(\theta)(M+m)}\,\omega o_I+\eta_e, \label{eq:oidotSegway}\\
f&=\frac{-\mathbb{I}(\theta)}{\cos\theta}\left(k_p\eta_e+k_d\omega+k_Io_I\right)-\frac{\mathbb{I}(\theta)}{\mathbb{I}_p}k_{cd}v, \label{eq:PIDSegway}
\end{align}
where $\mathbb{I}(\theta)\eta_e=\sin(\theta+\beta)$ stabilizes the vertical upright configuration $(\theta+\beta)=0$ semi-globally and exponentially.

\subsubsection{Simulation Verification}
We simulate performance of the controller (\ref{eq:oidotSegway}) -- (\ref{eq:PIDSegway}) for vertical stabilization of an inverted pendulum on a cart, which is an approximation of a segway type balancing vehicle. Nominal parameters used for the controller were $M=6.5\, kg$, $m=0.5\, kg$, $L=0.3\, m$ and $\mathbb{I}_p=0.09\, kgm^2$. To demonstrate the robustness of the controller system parameters were chosen to be $50\%$ different from nominal parameters used for the controller. Angle of the surface is considered to be $\beta=30^0$.    

\begin{figure}[h!]
	\centering
	\begin{tabular}{ccc}
		\includegraphics[width=.3\textwidth]{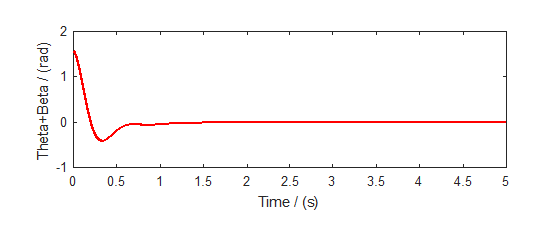}& \includegraphics[width=0.3\textwidth]{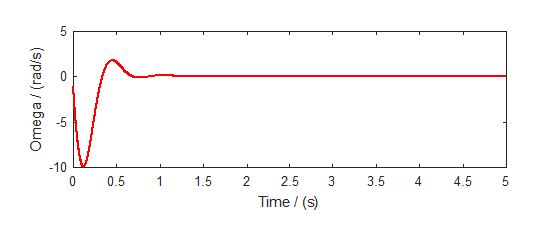} &
		\includegraphics[width=0.3\textwidth]{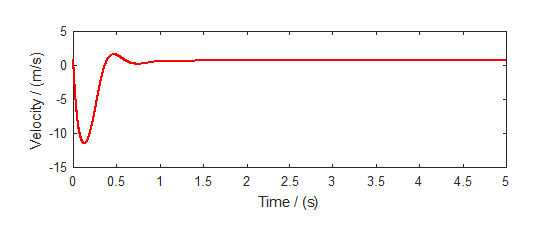}
	\end{tabular}
	\caption{The tilt angle $(\theta+\beta)$, the angular velocity $\dot{\theta}$, and the linear velocity $\dot{x}$ of the Inverted Pendulum on a Cart with the PID controller (\ref{eq:oidotSegway}) -- (\ref{eq:PIDSegway}) in the presence of parameter uncertainties as large as $50\%$. \label{Fig:IPC}}
\end{figure}

In figure--\ref{Fig:IPC} we present simulation results of tilt angle, tilt angular velocity and linear velocity in the presence of parameter uncertainties as large as $50\%$ and constant disturbances. Initial tilt angle, initial linear velocity and initial tilt angular velocity were considered to be $\theta(0)=60^0$, $v(0)=1\, ms^{-1}$ and $\omega(0)=1\, rads^{-1}$ respectively. The controller gains were chosen to be $k_p=10, k_d=6, k_I=2, k_{cd}=5$ and $k_{cp}=8$.

\subsection{Position Tracking of a Rolling Hoop on an Inclined Plane}\label{Secn:RollingHoop}

\begin{figure}[h!]
	\centering
	\begin{tabular}{c}
		\includegraphics[width=0.5\textwidth]{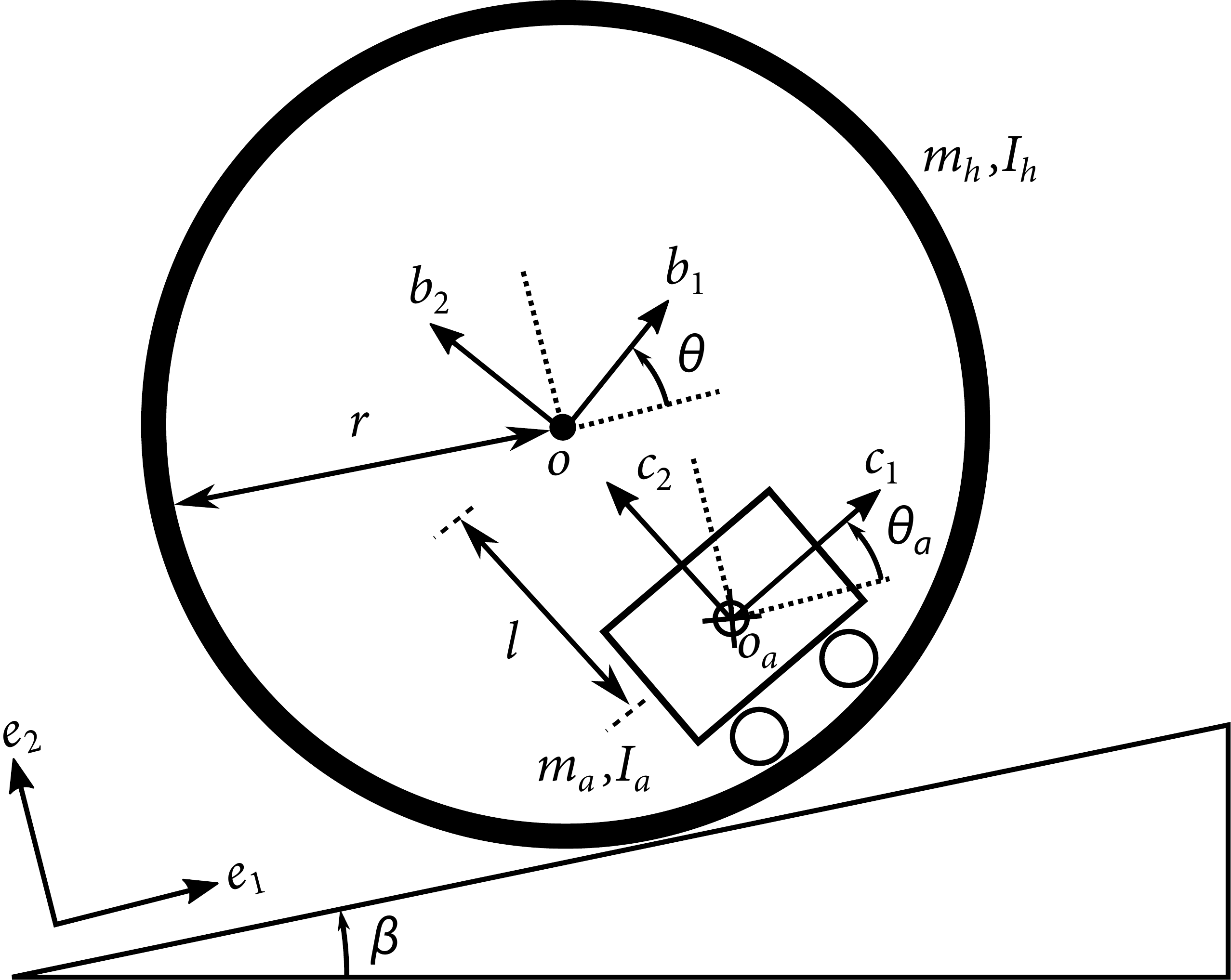}
	\end{tabular}
	\caption{Reference frames for geometric hoop robot analysis. The frames are well-defined for any actuation mechanism satisfying the constant-distance assumption.\label{Fig:HoopRefFrames}}
\end{figure}

Figure-\ref{Fig:HoopRefFrames} shows a hoop of radius $r$, rolling without slip on an inclined plane of constant inclination $\beta$ with respect to the horizontal plane. Let $\mathbf{e}=[\mathbf{e}_1\:\:\:\mathbf{e}_2]$ be an earth-fixed inertial frame with second axis pointing along the outward normal of the surface. Let $\mathbf{b}=[\mathbf{b}_1\:\:\:\mathbf{b}_2]$ be a reference frame fixed on the rolling hoop, with origin coinciding with the geometric center of the hoop. Let $\theta$ be the angle of rotation of the frame $\mathbf{b}$ with respect to the frame $\mathbf{e}$, and let $\omega=\dot{\theta}$. Let $m_h$ be the mass of the hoop and $\mathbb{I}_h$ be the inertia of the hoop. We assume that a mechanism of mass $m_a$ and moment of inertia $\mathbb{I}_a$ actuates the hoop. Let $\mathbf{c}=[\mathbf{c}_1\:\:\:\mathbf{c}_2]$ be a frame fixed to the actuation mechanism, with origin $o_a$ coinciding with the actuation center of mass and axis $\mathbf{c}_2$ pointing towards the center of the hoop. Let $\theta_a$ be the rotation angle of $\mathbf{c}$ with respect to $\mathbf{e}$ and $\omega_a=\dot{\theta}_a$.

Let $l$ be the distance from the geometric center of the hoop to the center of mass of the actuation mechanism. We restrict our attention to cases where $l$ remains constant. Under this assumption, the actuation mechanism evolves on the circle $\mathbb{S}$, and is completely characterized by configuration variable $\theta_a$. While this assumption is not strictly necessary to the rest of the analysis, it allows a compact characterization of the result.

We assume that the external forces and external moments acting on the hoop are due to the effect of gravity, reactions that arise as a consequence of the interaction with the actuation mechanism and the no-slip rolling constraint. The gravity force acting on the hoop is $f_g=-m_hg\,e_g$ where $e_g$ is the unit vector in the vertical direction, which can be written as $\sin\beta\,  \mathbf{e}_1 + \cos\beta\,  \mathbf{e}_2$ with respect to $\mathbf{e}$. Let $f_{\lambda}=[f_{{\lambda}_1} \:\:\: f_{{\lambda}_2}]^T \in \mathbb{R}^2$ be the $\mathbf{e}$-frame representation of the force that ensures the no-slip rolling constraints. The resultant moment acting on the hoop due to these constraint forces is $\tau_{\lambda}=rf_{{\lambda}_1}$. Let $f_c=[f_{c_1} \:\:\: f_{c_2}]^T \in \mathbb{R}^2$ be the $\mathbf{e}$-frame representation of the force acting on the hoop due to the interaction with the actuation mechanism and let $\tau_{f_c}$ be the resultant moment acting on the hoop. All moments are assumed to be taken with respect to the geometric center of the hoop. By Newton's third law, $-f_c$ and $-\tau_{f_c}$ are the reaction forces and moments acting on the actuation mechanism. For cart-type actuation mechanism the interaction between the hoop and the cart occurs only through no-slip constraints at the cart wheels, and hence $\tau_c\equiv 0$. For pendulum-type actuation, $\tau_{f_c}\equiv 0$.

Defining
{
	\begin{align*}
	M&\triangleq m_h+m_a , \\
	\tau^\omega_g&\triangleq rMg\sin{\beta}-\frac{m_a^2rl^2g}{\mathbb{I}_a+m_al^2}\cos{\theta_a}\sin{(\theta_a+\beta)},\\
	\tau^{\omega_a}_g&\triangleq \frac{m_arl\cos{\theta_a}}{\mathbb{I}_a+m_al^2}\tau^\omega_g-\left(\mathbb{I}_h+Mr^2-\frac{m_a^2r^2l^2}{\mathbb{I}_a+m_al^2}\,\cos^2\theta_a\right)\left(\frac{m_agl\sin({\theta_a+\beta})}{\mathbb{I}_a+m_al^2}\right) , \\
	B(\theta_a)&\triangleq \left(\frac{m_arl\cos{\theta_a}}{\mathbb{I}_a+m_al^2}-\frac{\left(\mathbb{I}_h+Mr^2-\frac{m_a^2r^2l^2}{\mathbb{I}_a+m_al^2}\,\cos^2\theta_a\right)}{(\mathbb{I}_a+m_al^2-m_arl\cos{\theta_a})}\right),
	\end{align*}} 
\noindent the expression (\ref{eq:ConstrainedMech}) that describe the motion of a non-holonomic mechanical system yields the following complete hoop robot equations of motion on the state space $\mathbb{S}\times\mathbb{R}\times\mathbb{R}\times\mathbb{S}\times\mathbb{R}$:
{
	\begin{align} 
	&\dot{\theta}=\omega,\\
	&\dot{o}=-r\omega,\\
	&\left(\mathbb{I}_h+Mr^2-\frac{m_a^2r^2l^2}{\mathbb{I}_a+m_al^2}\cos^2{\theta_a}\right)\dot{\omega}=-m_arl\sin{\theta_a}\omega_{a}^2+\tau^\omega_g+\tau_u,\\
	&\dot{\theta}_a={\omega}_a,\\
	&\left(\mathbb{I}_h+Mr^2-\frac{m_a^2r^2l^2}{\mathbb{I}_a+m_al^2}\,\cos^2\theta_a\right)\dot{\omega}_a= 
	-\frac{m^2_ar^2l^2\sin{\theta_a}\cos{\theta_a}}{\mathbb{I}_a+m_al^2}\omega_{a}^2+\tau^{\omega_a}_g+B(\theta_a)\tau_u .  \label{eq:CartDy}
	\end{align}}
The single control input, which appears in both the $\omega$ and $\omega_a$ equations, is defined as
\begin{align*}
\tau_u\triangleq \frac{(\mathbb{I}_a+m_al^2-m_arl\cos{\theta_a})}{\mathbb{I}_a+m_al^2}\left(\tau_c+\tau_{f_c}\right) .
\end{align*} 
We designate the output of the system to be the position of the hoop center, $o$.

The control task that we consider is to ensure that the output satisfies $\lim_{t \to \infty} o(t)=o_{\mathrm{ref}}(t)$ where $o_{\mathrm{ref}}(t)\in \mathbb{R}$ is a twice differentiable reference  and $\omega_{\mathrm{ref}}=-\dot{o}_{\mathrm{ref}}/r$. Let $o_e\triangleq (o-o_{\mathrm{ref}})$, $\omega_e\triangleq (\omega-\omega_{\mathrm{ref}})$. Note that in the special case of stabilizing the hoop at a point the references are constant: $o_{\mathrm{ref}}(t)\equiv \mathrm{const}$ and $\omega_{\mathrm{ref}}\equiv 0$.

Differentiating these quantities one sees that the error dynamics of the system take the explicit form
\begin{align}
\dot{o}_e&=-r\omega_e,\label{eq:Erroro}\\
\mathbb{I}({\theta_a})\dot{\omega}_e&=-m_arl\sin{\theta_a}\omega_{a}^2+\tau_g^{\omega}-\tau_{\mathrm{ref}}+\tau_u,\label{eq:Erroromega}
\end{align}
where $\mathbb{I}(\theta_a)\triangleq \mathbb{I}_h+Mr^2-\frac{m_a^2r^2l^2}{\mathbb{I}_a+m_al^2}\,\cos^2\theta_a$, and
$\tau_{\mathrm{ref}}\triangleq\mathbb{I}(\theta_a)\dot{\omega}_r$.
We notice that the error dynamics evolve on the tangent bundle to the circle, $T\mathbb{S}$, with output $y=o_e$ evolving on the Lie-group $\mathbb{R}$.
The natural notion of differentiation on the tangent space of a Riemannian manifold is the Levi-Civita connection, $\nabla$. As discussed in Section-\ref{Secn:MechSysOnS}, the unique Levi-Civita connection on $\mathbb{S}$ corresponding to the kinetic energy induced by the inertia tensor $\mathbb{I}$ is explicitly given by
\begin{align*}
\mathbb{I}({\theta_a})\nabla_{\zeta}\eta=\mathbb{I}({\theta_a})d\eta(\zeta)+\frac{m_a^2r^2l^2\sin{(2\theta_a)}}{2(\mathbb{I}_a+m_al^2)}\zeta\eta,
\end{align*}
for any $\zeta,\eta\in T_{g_a}\mathbb{S}$. Setting $\zeta$ to $\omega_a$ and $\eta$ to $\omega_e$ gives
\begin{align}
\mathbb{I}({\theta_a})\nabla_{\omega_a}\omega_e = \mathbb{I}({\theta_a})\dot{\omega}_e + \frac{m_a^2r^2l^2\sin{(2\theta_a)}}{2(\mathbb{I}_a+m_al^2)}\omega_a \omega_e.\label{eq:LeviCivitaS2}
\end{align}
Comparing with (\ref{eq:Erroromega}) and (\ref{eq:LeviCivitaS2})  it is clear that since the term
\begin{align*}
\frac{m_a^2r^2l^2\sin{(2\theta_a)}}{2(\mathbb{I}_a+m_al^2)}\omega_a \omega_e
\end{align*}
is absent (\ref{eq:Erroromega}) does not represent a simple mechanical system.
This prevents the straightforward use of the nonlinear PID controller proposed in Section-\ref{Secn:NLPID}.

However we notice that if we choose the regularizing plus potential shaping control
\begin{align}
\tau_u&=-\frac{m_a^2r^2l^2\,\sin{2\theta_a}}{2(\mathbb{I}_a+m_al^2)}\,\omega_a\omega_e+\frac{m_a^2rl^2g\sin{(2\theta_a)}}{2(\mathbb{I}_a+m_al^2)}+\tilde{\tau}_u,\label{eq:TransformingControl}
\end{align}
\noindent then the error dynamics of the system (\ref{eq:Erroromega}) can be re-written as,
\begin{align}
\mathbb{I}({\theta_a})\nabla_{\omega_a}\omega_e&=\Delta_h+\tau_{h}(\omega_e,\omega_a)+\tilde{\tau}_u.\label{eq:IntrinsicOmegaeDy}
\end{align}
The second term in the above equation (\ref{eq:TransformingControl}) shapes the potential energy of the error dynamics.
Here $\Delta_h$ represents the effects due to the ignorance of the inclination of the rolling surface and the omission of the term $\tau_{\mathrm{ref}}$ in the controller while
\begin{align*}
\tau_{h}(\omega_e,\omega_a)&\triangleq-m_arl\sin{\theta_a}\omega_{a}^2.
\end{align*}
In similar fashion we find that the actuation system dynamics can also be expressed as
{\small
	\begin{align}
	\mathbb{I}\nabla_{\omega_a}\omega_a&=\tilde{\tau}^{\omega_a}_g+\Delta_a+\tau_{a}(\omega_e,\omega_a)+B(\theta_a)\tilde{\tau}_u.\label{Eq:InnerCart}
	\end{align}}
\noindent where $\Delta_a$ represents modeling errors and disturbances. Here 
\begin{align*}
\tilde{\tau}^{\omega_a}_g&\triangleq {\tau}^{\omega_a}_g +B(\theta_a)\left(\frac{m_a^2l^2g\sin{(2\theta_a)}}{2(\mathbb{I}_a+m_al^2)}\right),\\
\tau_{a}(\omega_e,\omega_a)&\triangleq -B(\theta_a)\left(\frac{m_a^2r^2l^2\,\sin{2\theta_a}}{2(\mathbb{I}_a+m_al^2)}\,\omega_a\omega_e\right).
\end{align*}
Notice that $\tau_h(\zeta,\eta)$ and $\tau_a(\zeta,\eta)$ are bilinear in the two velocity arguments $\zeta,\eta\in\mathbb{R}$. 
When $\omega_e\equiv 0$ it can be shown that there exists an equilibrium for the actuation mechanism (\ref{Eq:InnerCart}) for any surface of inclination $\beta \in (-\pi/2,\pi/2)$ if the system parameters satisfy, $\frac{m_al}{Mr}\geq \sin\beta$. Without loss of generality we assume that the actuation mechanism is chosen such that this condition is satisfied for a certain operating region of $\beta$.

Note that the combination of the error dynamics of the system (\ref{eq:IntrinsicOmegaeDy}) and the actuation mechanism dynamics (\ref{Eq:InnerCart}) can be considered as an interconnected under actuated mechanical system where each sub system is fully actuated with respect to the common input $\tilde{\tau}_u$. The system evolves on $(\mathbb{S}\times \mathbb{R})\times (\mathbb{S}\times\mathbb{R})$ with the output $y=o_e$ evolving on the Lie-group $\mathbb{R}$.
For this interconnected mechanical system the nonlinear PID controller (\ref{eq:GeneralI}) -- (\ref{eq:GeneralPID}) reduces to,
\begin{align}
&\mathbb{I}({\theta_a})\nabla_{\omega_a}o_I=\mathbb{I}({\theta_a})\eta_e,\label{eq:IntegratorHoop}\\
&\tilde{\tau}_u=-\mathbb{I}({\theta_a})(k_p\eta_e+k_d\omega_e+k_Io_I+k_cB^{-1}(\theta_a)\omega_a),\label{eq:PIDHoop}
\end{align}
where, $\eta_e=-o_e$.

\subsubsection{Simulation Verification}
Here we present simulation results that demonstrate the effectiveness of the nonlinear PID tracking controller (\ref{eq:IntegratorHoop}) -- (\ref{eq:PIDHoop}) for a rolling hoop on an inclined plane. The actuation mechanism we consider is of the type of a cart or a pendulum. In order to demonstrate the robustness of the controller we choose the system parameters in the simulations to be $50\%$ different from the nominal parameters used in the controller. 

The nominal parameters for the hoop were chosen to be $m_h=1.00\, kg$ and $\mathbb{I}_h=0.021\, kgm^2$ while the outer radius of the hoop was chosen to be $r=0.18\, m$. These parameters were chosen to correspond to a hoop made of plastic (density $850\, Kgm^3$) with thickness of $3\, mm$. The nominal parameters for the actuation appendage was chosen to be $m_a =3.28\, kg$ and $\mathbb{I}_a=0.035\, kgm^2$. The distance from the geometric center of the hoop to the center of mass of the actuation mechanism was chosen to be $l=0.14\, m$. For these parameters we find that the maximum inclination of the rolling plane for which an equilibrium exists for the actuation mechanism is $\beta_{\mathrm{max}}=36^0$.
Thus in the simulations the hoop is assumed to roll on a inclined plane of $20^0< \beta_{\mathrm{max}}$. We stress that the implementation of the controller (\ref{eq:IntegratorHoop})--(\ref{eq:PIDHoop}) does not require the knowledge of the angle of inclination of the rolling surface.

\begin{figure}[h!]
	\centering
	\begin{tabular}{ccc}
		\includegraphics[width=0.32\textwidth]{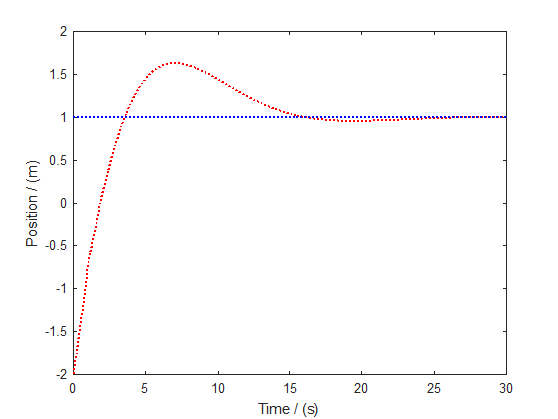}& \includegraphics[width=0.32\textwidth]{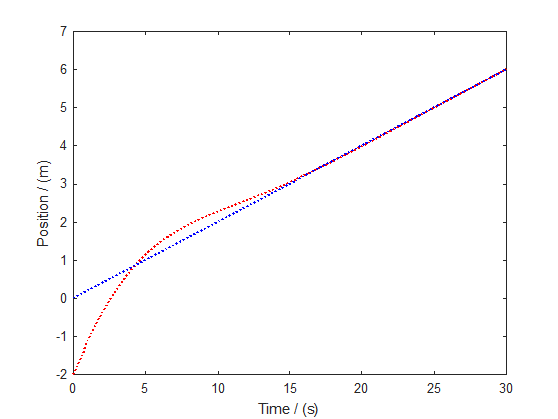} & \includegraphics[width=0.32\textwidth]{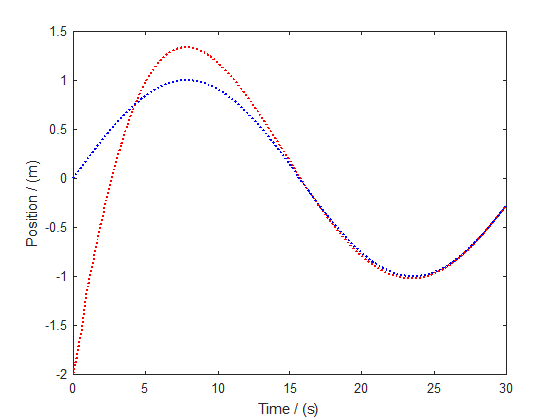}\\
		\includegraphics[width=0.32\textwidth]{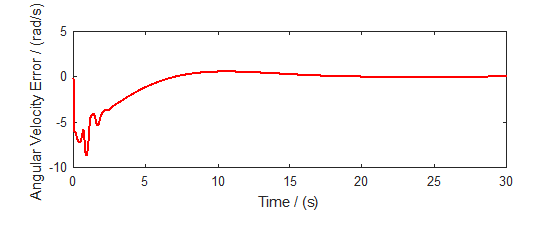}& \includegraphics[width=0.32\textwidth]{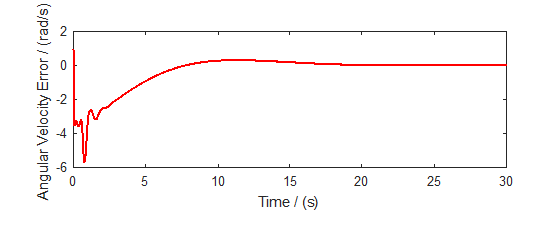}&
		\includegraphics[width=0.32\textwidth]{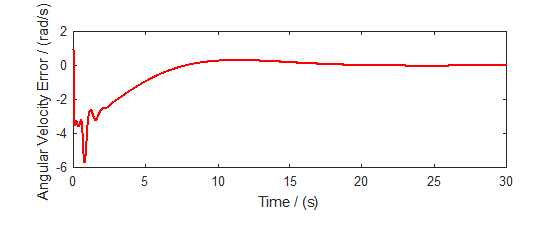}\\
		\includegraphics[width=0.32\textwidth]{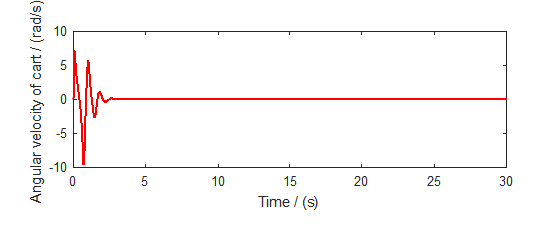}& \includegraphics[width=0.32\textwidth]{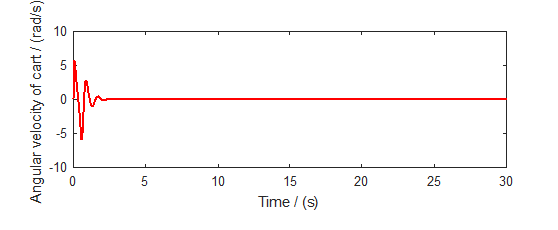}&
		\includegraphics[width=0.32\textwidth]{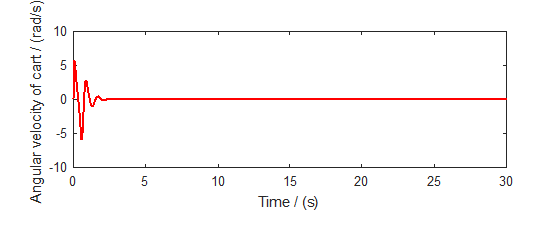}\\(a) Fixed point & (b) Linear velocity profile & (c) Sinusoidal velocity profile
	\end{tabular}
	\caption{The position followed by the geometric center of the hoop, $o(t)$, the angular velocity tracking error, $\omega_e(t)$, and the angular velocity of the actuation mechanism, $\omega_a(t)$  for the PID controller (\ref{eq:IntegratorHoop}) -- (\ref{eq:PIDHoop}) in the presence of parameter uncertainties as large as $50\%$. The blue curve in the first row of figures shows the reference position while the red curve shows the actual position of the center of the hoop.\label{Fig:RollingHoopPosition}}
\end{figure}

We present results for: a.) stabilizing the hoop at a point, b:) tracking a linear velocity, and c.) tracking a sinusoidal velocity. The simulation results are presented in Figure-\ref{Fig:RollingHoopPosition}. Figure-\ref{Fig:RollingHoopPosition}(a) -- figure-\ref{Fig:RollingHoopPosition}(c) shows the reference position and actual position of the geometric center of the hoop, spatial angular velocity error of the hoop and the spatial angular velocity of the actuation mechanism for stabilizing at a point, tracking a linear velocity and tracking a sinusoidal velocity respectively.
In all simulations the initial position of the hoop was assumed to be at $o(0)=[-2\:\:\:r]^T\,m$ and the initial spatial angular velocity for the hoop and inner cart were chosen as $\omega(0)=-0.1\,rads^{-1}$ and $\omega_a(0)=0.1\,rads^{-1}$ respectively. In all simulations the controller gains were chosen to be $k_p=16, k_d=7, k_I=4, k_c=0.1$.

\subsection{Trajectory Tracking of a Rolling Sphere on an Inclined Plane}\label{Secn:RollingSphere}
Consider an inertially asymmetric balanced sphere rolling on a surface of unknown constant inclination. The motion of the sphere is generated by an actuation mechanism. 
We model the system as a composite system composed of the rigid sphere plus rigid actuation mechanism. Forces and moments acting at the point of contact of the actuation rigid body with the sphere couples the dynamics of the actuation mechanism to the dynamics of the sphere. 
We will use intrinsic Newton's equations to model each of the rigid bodies and invoke Newton's $3^{\mathrm{rd}}$ law to couple the forces and moments acting at the contact point.
Specifically we will derive the equations of motion of the system by considering each body separately using (\ref{eq:ConstrainedMech}).

The configuration space of the rolling sphere is $\mathbb{R}^2\times SO(3)\subset SE(3)$. Let $\mathbf{e}$ be an orthonormal inertial frame with the third axis coinciding with the outward normal to the surface of rolling and let $\mathbf{b}$ be an orthonormal frame fixed to the sphere with origin coinciding with the geometric center of the sphere. 
Let $o(t)$ be the representation of the geometric center of the sphere with respect to the frame $\mathbf{e}$ and let $\mathbf{b}=\mathbf{e}R(t)$ where $R(t)\in SO(3)$. Let $m_b$ be the total mass of the spherical shell and let $\mathbb{I}_b$ be the inertia tensor of the spherical shell with respect to $\mathbf{b}$. We will not assume that the metric induced by $\mathbb{I}_b$ is bi-invariant. That is, we will not assume $\mathbb{I}_b$ to be of the form $\mu I_{3\times 3}$ for some positive $\mu$.

We assume that there are no external forces and moments acting on the sphere except for gravity and the constraints. Similarly we will assume that actuation mechanism only feels gravity, the reaction forces, $-f_{c}$, reaction moments, $-\tau_{c}$, and the moments due to the reaction forces, $-\tau_{f_{c}}$, that arise due to its interaction with the sphere. All moments are taken with respect to the geometric center of the sphere. 
The gravity force acting on the sphere is $f_g=-m_bg\,e_g$ where $-e_g$ is the unit direction that coincides with the direction of gravity in the inertial frame $\mathbf{e}$.
Let 
$f_c$ be the resultant of the reaction forces, let $\tau_c$ be the resultant of the reaction moments, and let $\tau_{f_c}$ be the resultant moment on the sphere due to the reaction forces on the sphere that arises as a consequence of the  interactions it has with the actuation mechanism. 
The constraint that the sphere is always in contact with the surface implies that $e_3\cdot \dot{o}=0$.
Let $p$ be an arbitrary point on the surface of the sphere with representation $x(t)$ and $X$ with respect to $\mathbf{e}$ and $\mathbf{b}$ respectively. Then since
$\dot{x}=\dot{o}+R({\Omega}\times X)=\dot{o}+({\omega}\times RX)$, the no-slip rolling implies $\dot{x}=0$ at the point of contact of the sphere with the surface on which it rolls and hence we have that the velocity constraints on the sphere are given by
\begin{align*}
\dot{o}+r{e}_3\times\omega&=0,\\
e_3\cdot \dot{o}&=0,
\end{align*}
in terms of the spatial velocities on $SO(3)$. 

Let $\bar{X}_i$ be the representation of the center of mass of the actuation mechanism with respect to the frame $\mathbf{b}$ that is fixed to the sphere. That is if $o_i$ is the representation of the center of mass of the actuation mechanism with respect to the earth fixed frame $\mathbf{e}$ then $o_i=o+R\bar{X}_i$ and $\dot{o}_i=-(r\hat{e}_3+\widehat{R\bar{X}_i}){\omega}+R\dot{\bar{X}}_i$. 
Let us introduce the following angular velocities: $\mathbf{b}_i=\mathbf{b}R_{b_i}=\mathbf{e}RR_{b_i}=\mathbf{e}R_i$ where $R_{b_i}\in SO(3)$. Let $\dot{R}_{i}=R_{i}\widehat{\Omega}_{i}=\widehat{\omega}_{i}R_{i}$, and $\dot{R}_{b_i}=R_{b_i}\widehat{\Omega}_{b_i}=\widehat{\omega}_{b_i}R_{b_i}$. Here $\Omega_i$ is the body velocity and $\omega_i$ is the corresponding spatial velocity of the  actuation mechanism while ${\Omega}_{b_i}$ is the body velocity with respect to the frame $\mathbf{b}$ fixed on the sphere. Note that these versions are related by ${\omega}_i={\omega}+R{\omega}_{b_i}={\omega}+R_i{\Omega}_{b_i}$.  
The convenient choice of variable to use to represent the actuation mechanism will depend on the type of mechanism being used.
Let $-f_{i}=-f^u_{i}-f^e_{i}$ be the interaction force and let $-\tau_{i}$ be the interaction moment  acting on the body due to its interaction at $P_{i}$.  Note that
the interaction force $f_{i}$ also gives rise to a moment on the actuation mechanism as well as the sphere. That is if $X_{i}$ is the representation of the point $P_{i}$ in the frame $\mathbf{b}_i$ it effects a moment $-R_iX_{i}\times f_{i}$ on the actuation mechanism about its center of mass and a moment $\tau_{f_{i}}=(R\bar{X}_i+R_iX_{i})\times f_{i}$ on the sphere about the geometric center of the sphere. 

Let the mass of the actuation mechanism be $m_i$ and its moment of inertia tensor in the frame $\mathbf{b}_i$ about its center of mass be $\mathbb{I}_i$. Here we consider that the distance between geometric center of the sphere and center of mass of the actuation mechanism to be a constant denoted as $l$. Then dynamic equation of the sphere can be expressed as,
\begin{align}
\mathbb{I}_b^R\nabla^b_{\omega}\omega+\mathbb{I}_s\dot{\omega}
+ \sum_{i=1}^n{e}_3\times\left(\mathbb{I}^{R_i}_{\tilde{V}_i}{({e}_3\times \dot{\omega})}\right)
&=B\left(\tau_{c}+\tau_{f_c}\right)+\tau_v+\tau_g, \label{eq:BallDynamicsVehicle}
\end{align}
where we have set 
\begin{align*}
\mathbb{I}_s&\triangleq -(m_b+m_I)r^2\widehat{e}_3^2\\
\mathbb{I}_{V_i}&\triangleq \left(\mathbb{I}_{i}-m_il^2\widehat{e}_3^2\right)\\
\mathbb{I}_{\tilde{V}_i}&\triangleq -r^2m_i^2l^2\widehat{e}_3\mathbb{I}_{V_i}^{-1}\widehat{e}_3 \\
\tau_v&\triangleq\sum_{i=1}^nm_irl\widehat{e}_3R_i\widehat{e}_3\mathbb{I}_V^{-1}R_i^T\widehat{\omega}_i\mathbb{I}_V^{R_i}\omega_i+m_irl\widehat{e}_3\widehat{\omega_i}^2R_ie_3\\
\tau_g&\triangleq -r(m_b+m_I)g\,e_3\times e_g+\frac{g}{r}\sum_{i=1}^n{e}_3\times\mathbb{I}_{\tilde{V}_i}^{R_i}e_g\\
B&\triangleq \sum_{i=1}^n\left(m_i rl\widehat{e}_3R_i\widehat{e}_3{\mathbb{I}_{V_i}}^{-1}R_i^T+I_{3\times 3}\right),
\end{align*}
and 
$\mathbb{I}_{\tilde{V}_i}^{R_i}\triangleq R_i\mathbb{I}_{\tilde{V}_i}R_i^T$.
We also find that the governing equation for the actuation mechanism is given by 
\begin{align}
\mathbb{I}_{V_i}^{R_i}\nabla^{V_i}_{{\omega_i}}{\omega_i}&=-m_irl\,R_ie_3\times{e}_3\times\dot{\omega}+m_igl\,R_ie_3\times e_g
-(\tau_{c}+\tau_{f_{c}}),\label{eq:ithbodyomegaVehicle}
\end{align}
where  $\nabla^{V_i}$ is the unique Levi-Civita connection corresponding to the left-invariant Riemannian metric induced by $\mathbb{I}_{V_i}$.

From the equation (\ref{eq:ConstrainedMech}) that describe the motion of constrained mechanical systems presented  section-\ref{Secn:CovDerivative} it follows that the dynamics of the rolling sphere is given by
{\begin{align*}
	\dot{R}&=\widehat{\omega}R,\\
	\dot{o}&=r\omega\times e_3,\\
	\mathbb{I}_\alpha^R\nabla^\alpha_{\omega}\omega+\left(\mathbb{I}_{s}
	+\sum_{i=1}^n\widehat{e}_3 \mathbb{I}^{R_i}_{\tilde{V}_i}\widehat{e}_3 \right)\dot{\omega}&=\tau_\alpha+\tau_g+\mathbb{I}\Delta_d+B\tau_u,
	\end{align*}
}
where $\mathbb{I}\Delta_d$ represent unmodelled moments that may arise for instance due to the uncertainties in the knowledge of the inclination of the rolling surface.

Notice that the second term on the left hand side of the last of the above equations prevents the system from being treated as a left-invariant mechanical system on $SO(3)$. 
However if we define the singular right-invariant Riemannian metric $\langle\langle\cdot,\cdot\rangle\rangle_{s}$ on $SO(3)$ by the relationship $\langle\langle \widehat{\zeta}R,\widehat{\eta}R\rangle\rangle_{s}=\zeta^T\mathbb{I}_{s}\eta$ then the sphere dynamics can be written down as
{
	\begin{align}
	&\mathbb{I}_s\nabla^s_{\omega}\omega
	+\mathbb{I}_\alpha^R\nabla^\alpha_{\omega}\omega+{e}_3\times \sum_{i=1}^n\left(\mathbb{I}^{R_i}_{\tilde{V}_i}\nabla^{\tilde{V}_i}_{ \omega_i}{({e}_3\times {\omega})}\right)
	=\nonumber\\
	&-\sum_{i=1}^n\frac{1}{2}e_3\times \left(\mathbb{I}^{R_i}_{\tilde{V}}({\omega_i}\times{({e}_3\times {\omega})})+\mathbb{I}^{R_i}_{\tilde{V}}{\omega_i}\times{({e}_3\times {\omega})}+\mathbb{I}^{R_i}_{\tilde{V}_i}{({e}_3\times {\omega})}\times{\omega_i}\right)
	+\mathbb{I}_s\omega\times \omega
	+\tau_\alpha+\tau_g+\mathbb{I}\Delta_d+B\tau_u, \label{eq:ConstraintBallDynamicsVehicle}
	\end{align}
}
where $\nabla^{\tilde{V}_i}$ is the unique Levi-Civita connection corresponding to the left invariant Riemannian metric induced by $\mathbb{I}_{\tilde{V}_i}$.

We consider the case where $\tau_u$ has three degrees of freedom.
The control problem that we solve in this section is that of ensuring $\lim_{t \to \infty} o(t)=o_{\mathrm{ref}}(t)$ where $o_{\mathrm{ref}}(t)\in \mathbb{R}^2\subset \mathbb{R}^3$ is a twice differentiable and satisfies $e_3\cdot o_{\mathrm{ref}}(t)\equiv r$ for all $t$ while ensuring $\lim_{t \to \infty} \omega(t)=\omega_{\mathrm{ref}}(t)$ where $\omega_{\mathrm{ref}}(t)$ is such that 
$e_3\cdot\omega_{\mathrm{ref}}\equiv 0$ and $\dot{o}_{\mathrm{ref}}=r\,\omega_{\mathrm{ref}}\times {e}_3$. 
Let $o_e\triangleq (o-o_{\mathrm{ref}})$, $\omega_e\triangleq (\omega-\omega_{\mathrm{ref}})$.
Then we see that the error dynamics can be expressed as
{
	\begin{align*}
	&\dot{o}_e=-r\,{e}_3\times\omega_e,\\
	&\mathbb{I}_{s}\nabla^{s}_{{\omega}_e}{\omega}_e+\mathbb{I}_\alpha^R\nabla^{\alpha} _{{\omega}_e}{\omega}_e+{e}_3\times\sum_{i=1}^n \left(\mathbb{I}^{R_i}_{\tilde{V}_i}\nabla^{\tilde{V}_i}_{\omega_i}{(\widehat{e}_3 {\omega_e})}\right)=\mathbb{I}_s\omega_e\times \omega_e
	+B\tau_u-\tau_{\mathrm{ref}}+\mathbb{I}\Delta_d+\tau_\alpha+\tau_g+\tau_{\tilde{V}} ,
	\end{align*}
}
where
{ \begin{align*}
	\tau_{\mathrm{ref}}&=
	\left(\mathbb{I}_\alpha^R+\mathbb{I}_{s}
	+\widehat{e}_3 \left(\sum_{i=1}^n\mathbb{I}^{R_i}_{\tilde{V}_i}\right)\widehat{e}_3 \right)\dot{\omega}_{\mathrm{ref}}
	-\mathbb{I}_\alpha^R{{\omega}_{\mathrm{ref}}}\times{{\omega}_e}-\mathbb{I}_\alpha^R{{\omega}_e}\times{{\omega}_{\mathrm{ref}}}-\mathbb{I}_\alpha^R{{\omega}_{\mathrm{ref}}}\times{{\omega}_{\mathrm{ref}}},
	\end{align*}}
and {
	\begin{align*}
	\tau_{\tilde{V}}=-\frac{1}{2}e_3\times \sum_{i=1}^n&\left(\mathbb{I}^{R_i}_{\tilde{V}_i}({\omega_i}\times{({e}_3\times {\omega_e})})
	+\mathbb{I}^{R_i}_{\tilde{V}_i}{\omega_i}\times{({e}_3\times {\omega_e})}+\mathbb{I}^{R_i}_{\tilde{V}_i}{({e}_3\times {\omega_e})}\times{\omega_i}\right).
	\end{align*}}
These define a dynamic system on $\mathbb{R}^2\times so(3)$ and is the natural extension of the error dynamics for position tracking for the split mechanical system (\ref{eq:ConstraintBallDynamicsVehicle}).

The basic idea behind the conventional linear PID controller is that it reduces the error dynamics of a tracking problem, for a double integrator, to the form $\dot{e}_I=e$, $\ddot{e}=-k_pe-k_d \dot{e}-k_Ie_I$. We have shown in \cite{MaithripalaAutomatica} that for a general mechanical system on a manifold, which is in essence an intrinsic double integrator, the appropriate PID controller is $\nabla^{\alpha} _{\dot{e}}e_I=\mathrm{grad}^{\alpha}  V(e)$, $u=-k_p\mathrm{grad}^{\alpha}  V(e)-k_d\dot{e}-k_Ie_I$ and that it results in the intrinsic error dynamics $\nabla^{\alpha} _{\dot{e}}e_I=\mathrm{grad}^{\alpha}  V(e)$, $\nabla^{\alpha} _{\dot{e}}\dot{e}=-k_p\mathrm{grad}^{\alpha}  V(e)-k_d\dot{e}-k_Ie_I$ for some polar Morse function $V(e)$ on the configuration space. Thus given the split mechanical structure of the rolling ball we are motivated to propose the intrinsic nonlinear PID controller
\begin{align}
\mathbb{I}_s\nabla^s_{\omega_e}{o}_I&+\mathbb{I}_\alpha^R\nabla^{\alpha} _{\omega_e}{o}_I+{e}_3\times\sum_{i=1}^n\left(\mathbb{I}_{\tilde{V}_i}^{R_i}\nabla^{\tilde{V}_i} _{\omega_i}({e}_3\times {o}_I)\right)=
\mathbb{I}_s{\eta^s_e}+\mathbb{I}_\alpha^R{\eta^\alpha_e}+{e}_3 \times\sum_{i=1}^n\left(\mathbb{I}_{\tilde{V}_i}^{R_i}{ \eta^{\tilde{V}_i}_e}\right),\label{eq:Integrator}
\end{align}
\begin{align}
\tau_u&=-B^{-1}\left(\tau_\alpha+\tau_{\tilde{V}}+\mathbb{I}_s\omega_e\times \omega_e -\tau_{\mathrm{ref}}
+\frac{g}{r}\sum_{i=1}^n{e}_3\times\mathbb{I}_{\tilde{V}_i}^{R_i}e_3
\right.\nonumber\\&\:\:\:\:
+\mathbb{I}^R_\alpha(k_p\eta^\alpha_e +k_d\omega_e+k_Io_I)
+\mathbb{I}_s(k_p\eta^s_e +k_d\omega_e+k_Io_I)
\nonumber\\&\:\:\:\:
\left.+\widehat{e}_3 \left(\sum_{i=1}^n\mathbb{I}^{R_i}_{\tilde{V}_i}(k_p\eta^{\tilde{V}_i}_e +k_d\widehat{e}_3\omega_e+k_I\widehat{e}_3o_I)
\right)\right)+\mathbb{I}^{R_i}_{{V}_i}k_{cd}\omega_i,\label{eq:PIDcontroller}
\end{align}
where $\eta^\nu_e  =\mathrm{grad}^{\nu}  V_{\nu} (o_e)$ for a quadratic (a polar morse) function on $\mathbb{R}^2\subset \mathbb{R}^3$. Let
$V_{\nu} (o_e)\triangleq \frac{1}{2r}o_e\cdot o_e$.  Then since
$\dot{V}_\nu(o_e)=\langle dV_\nu,\omega_e\rangle=\langle\langle \eta_e^\alpha,\omega_e\rangle\rangle_\nu=\langle \mathbb{I}_\alpha^R\eta_e^\alpha,\omega_e\rangle$ the gradient of $V_\nu$, denoted by $\eta^\nu_e$, is explicitly given by
$\mathbb{I}^R_\nu\eta^\nu_e=dV_\nu =\widehat{e}_3 o_e\in \mathbb{R}^2\subset\mathbb{R}^3\simeq so(3)$ for $\nu=\alpha$ and $\nu=s$. For $\nu=\tilde{V}_i$ we let $V_{\tilde{V}_i} (o_e)\triangleq \frac{1}{2r}(\widehat{e}_3o_e)\cdot (\widehat{e}_3o_e)$ and since
$\dot{V}_{\tilde{V}_i}(o_e)=\langle dV_{\tilde{V}_i},\widehat{e}_3\omega_e\rangle=\langle\langle \eta_e^{\tilde{V}_i},\widehat{e}_3\omega_e\rangle\rangle_{\tilde{V}_i}=\langle \mathbb{I}_{\tilde{V}_i}^{R_i}\eta_e^{\tilde{V}_i},\widehat{e}_3\omega_e\rangle$ the gradient of $V_{\tilde{V}_i}$, denoted by $\eta^{\tilde{V}_i}_e$, is explicitly given by
$\mathbb{I}^{R_i}_{\tilde{V}_i}\eta^{\tilde{V}_i}_e=dV_{\tilde{V}_i} = o_e\in \mathbb{R}^2\subset\mathbb{R}^3\simeq so(3)$.

With this controller we see that the error dynamics are given by (\ref{eq:Integrator}) and
{
	\begin{align}
	&\dot{o}_e=-r\,{e}_3\times \omega_e\label{eq:Erroroe}\\
	&\mathbb{I}_{s}\nabla^{s}_{{\omega}_e}{\omega}_e+\mathbb{I}_\alpha^R\nabla^{\alpha} _{{\omega}_e}{\omega}_e+{e}_3\times \sum_{i=1}^n\left(\mathbb{I}^{R_i}_{\tilde{V}_i}\nabla^{\tilde{V}_i}_{ \omega_i}{(\widehat{e}_3 {\omega_e})}\right)=\mathbb{I}\Delta_d+\mathbb{I}\Delta_\epsilon \nonumber\\
	&\:\:\:\:-\mathbb{I}^R_\alpha\left(k_p\eta_e^{\alpha} +k_d\omega_e+k_Io_I\right)-\mathbb{I}_s\left(k_p\eta_e^s +k_d\omega_e+k_Io_I\right)
	-\widehat{e}_3\left(\sum_{i=1}^n\mathbb{I}^{R_i}_{\tilde{V}_i}(k_p\eta_e^{\tilde{V}_i} +k_d\widehat{e}_3\omega_e+k_I\widehat{e}_3o_I)\right)+B\mathbb{I}^{R_i}_{{V}_i}k_{cd}\omega_i,\label{eq:Errorwe}
	\end{align}
}
where $\mathbb{I}\Delta_\epsilon$ represents errors that arise due to the uncertainty in the inertia parameters. The error dynamics (\ref{eq:Integrator}) and (\ref{eq:Erroroe}) -- (\ref{eq:Errorwe}) evolve  on $\mathbb{R}^2\times so(3)\times so(3)$ and is the natural notion of error dynamics for a system that exhibits a split mechanical structure such as (\ref{eq:ConstraintBallDynamicsVehicle}). 

\subsubsection{Simulation Verification}
In this section we simulate the performance of the intrinsic nonlinear PID controller (\ref{eq:Integrator}) -- (\ref{eq:PIDcontroller}) for trajectory tracking for a sphere rolling on an inclined plane with unknown constant inclination.

For simulations the nominal mass of the spherical shell was chosen to be  $m=1.00\, kg$, the nominal inertia tensor of the spherical shell was chosen to be $\mathbb{I}_b=\mathrm{diag}\{0.0213, 0.0205, 0.0228\}\, kgm^2$, while the radius of the spherical shell was chosen to be $r=0.18\,m$. 
These parameters were chosen to correspond to a $3\, mm$ thickness shell made of plastic (density $850\, Kgm^{-3}$). 
To demonstrate robustness of the controller the system parameters used in the simulations were chosen to be $50\%$ different from the nominal parameters used for the controller. In all simulations the initial position of the sphere was assumed to be $o(0)=[2\:\:\:-2\:\:\:r]^T\,m$ and the initial angular velocity of the sphere was chosen to be $\omega(0)=[-0.1\:\:\:-0.2\:\:\:0.5]^T\,rads^{-1}$.

Consider a omni-directional wheel driven cart actuated  sphere. The total mass of the cart was chosen to be $m_i=3.28\,kg$, while its inertia tensor was chosen to be
\begin{align*}
\mathbb{I}_i&=\mathrm{diag}\{0.0353, 0.0378, 0.0368\}\,kgm^2.
\end{align*}
For these parameters one finds that the maximum inclination for which an equilibrium for the controlled cart exists is $25^0$.
The sphere is assumed roll on a $20^0$ inclined plane in the $y$-direction and is assumed to be unknown. A nominal value of $30^0$ inclination is assumed in the controller.

\begin{figure}[h!]
	\centering
	\begin{tabular}{ccc}
		\includegraphics[width=0.32\textwidth]{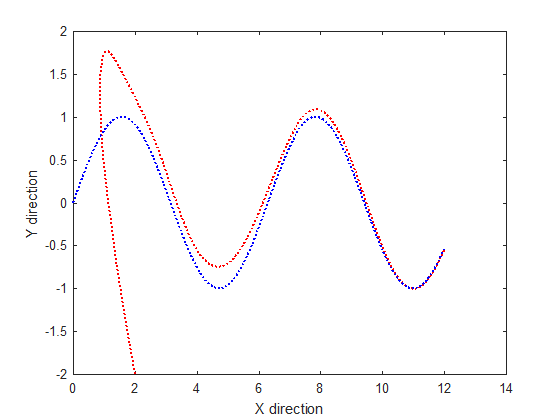}& \includegraphics[width=0.32\textwidth]{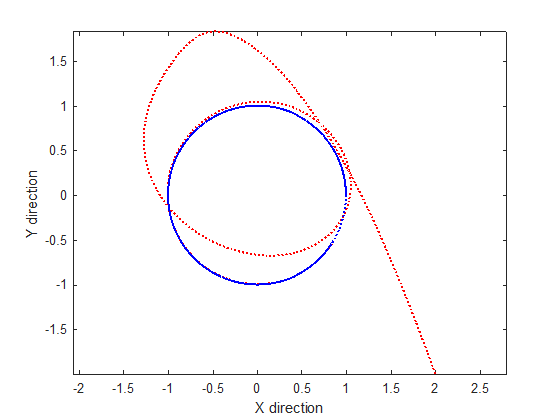}&
		\includegraphics[width=0.32\textwidth]{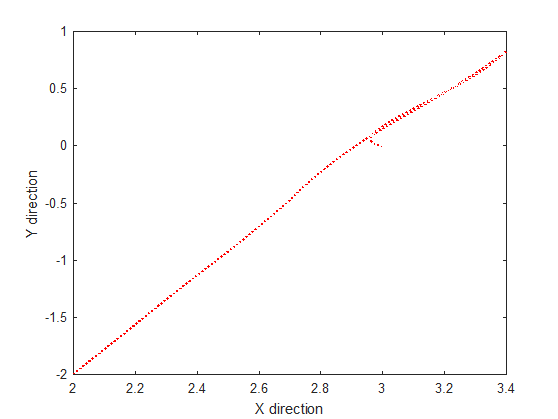}\\
		\includegraphics[width=0.32\textwidth]{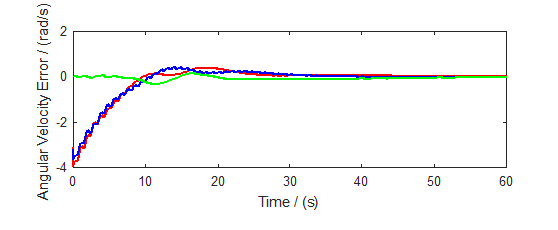}& \includegraphics[width=0.32\textwidth]{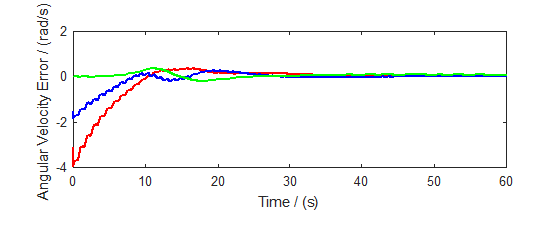}&
		\includegraphics[width=0.32\textwidth]{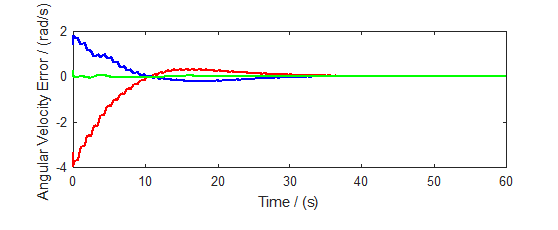}\\
		\includegraphics[width=0.32\textwidth]{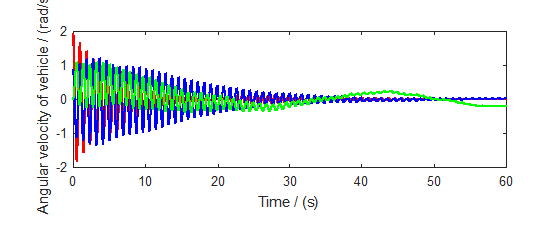}& \includegraphics[width=0.32\textwidth]{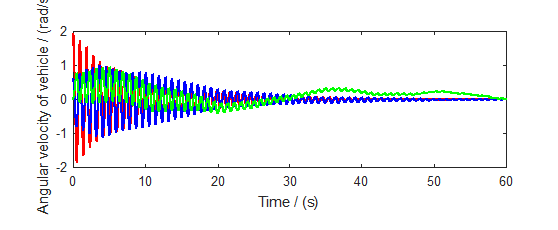}&
		\includegraphics[width=0.32\textwidth]{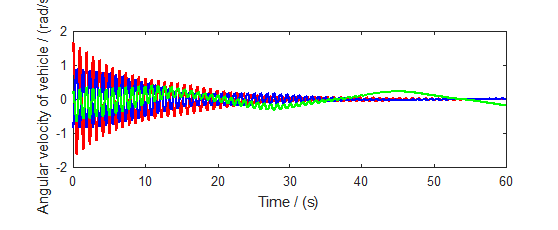}\\
		(a) Sinusoidal path & (b) Circular path &(c) Fixed point
	\end{tabular}
	\caption{The path followed by the center of mass $o(t)$ of the sphere, spatial angular velocity error, $\omega_e(t)$, and the spatial angular velocities of the cart, $\omega_{i}(t)$,  of the cart actuated sphere for the PID controller (\ref{eq:Integrator}) -- (\ref{eq:PIDcontroller}) in the presence of parameter uncertainties as large as $50\%$. The blue curve in the first row of figures shows the reference trajectory while the red curve shows the trajectory of the center of mass of the sphere.\label{Fig:RollingVehiclePath}}
\end{figure}

Simulation results are presented in figure-\ref{Fig:RollingVehiclePath} for tracking a sinusoidal path, circular path and a fixed point at $(3,0)$.
In all simulations the initial conditions used for the cart were $\omega_i(0)=[0.2\:\:\:-0.1\:\:\:0.1]^T\,rads^{-1}$.
The controller gains were chosen to be $k_p=100, k_d=60, k_I=10$.  

\subsection{Vertical Stabilization of a Spherical Pendulum}\label{Secn:SphericalPendulum}
\begin{figure}[h!]
	\centering
	\begin{tabular}{c}
		\includegraphics[width=0.4\textwidth]{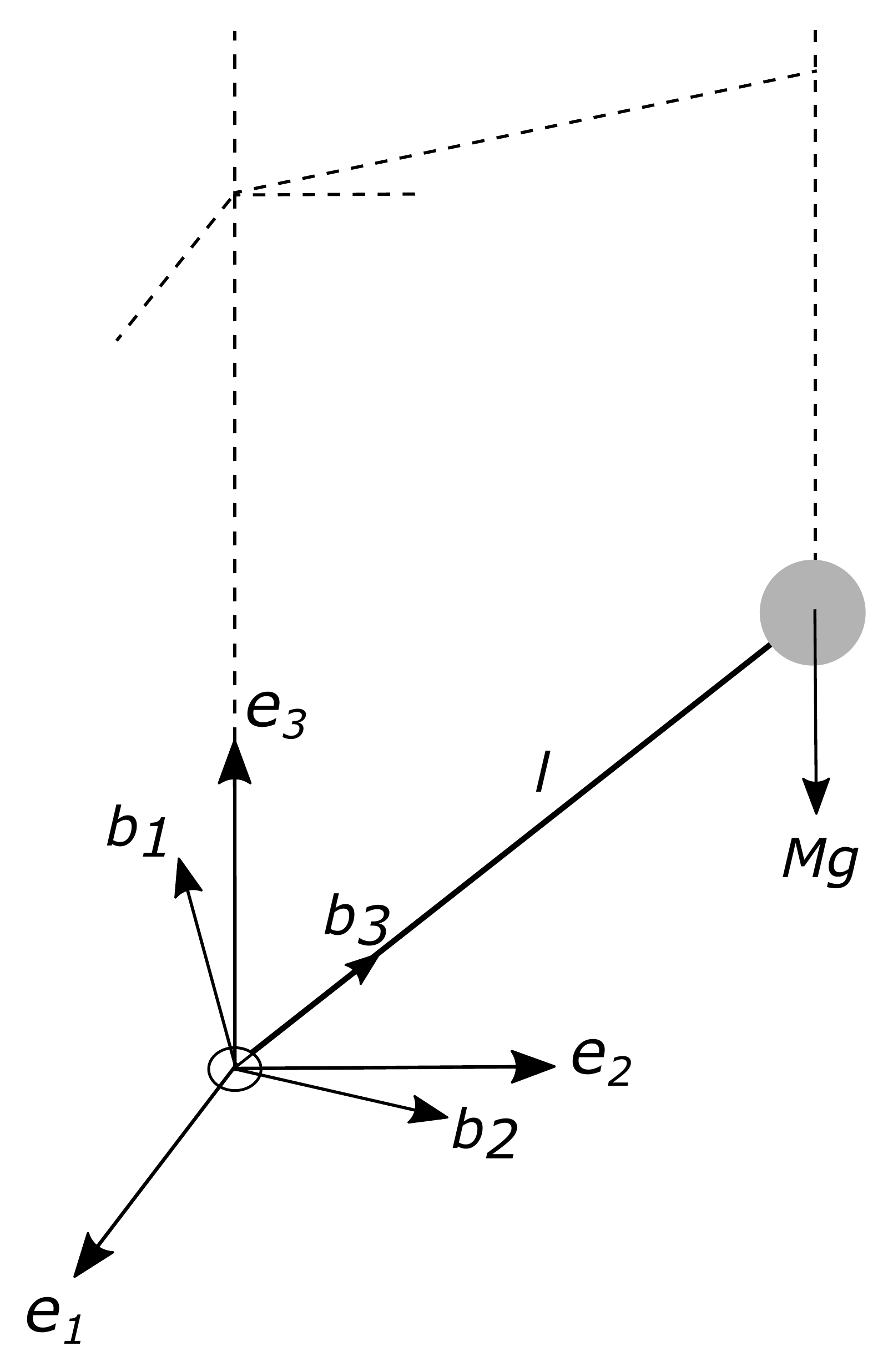}
	\end{tabular}
	\caption{Spherical Pendulum.\label{Fig:SphericalPendulum}}
\end{figure}
Consider the spherical pendulum shown in Figure-{\ref{Fig:SphericalPendulum}} that consists of a non-trivial spherical mass and a massless rod. Let $\mathbf{e}$ be an earth-fixed frame with origin $o_b$ coinciding with the pivot point of the pendulum. The $3^{\mathrm{rd}}$ axis is pointing vertically upward. Let $\mathbf{b}$ be a reference frame fixed to the spherical pendulum with origin coinciding with the origin $o_b$ of the earth fixed frame and let $\mathbf{b}=\mathbf{e}R(t)$ where $R(t)\in SO(3)$. Let $M$ be the mass of the pendulum, $l$ be the length of the rod, and $\mathbb{I}=\mathrm{diag}\{\mathbb{I}_1,\mathbb{I}_2,\mathbb{I}_3\}$ be the inertia tensor of the pendulum with respect to the $\mathbf{b}$ frame.

Let $\widehat{\Omega}=R^T\dot{R}$ be the skew-symmetric version of the body angular velocity of the pendulum and $\tau_u=[{\tau_u}_1\:\:\:\:{\tau_u}_2\:\:\:\:0]^T$ be the control moments.
The constraint that ensures that the pendulum does not rotate about the $\mathbf{b}_3$ axis, can be expressed as, 
\begin{align}
e_3\cdot{\Omega}=\Omega_3&=0. \label{eq:SPConstraint}
\end{align}
Thus the constraint distribution is defined by,
\begin{align*}
\mathcal{D}=\mathrm{Im}\left(\widehat{e}_3^2\right).
\end{align*}
and the dual projection maps onto the co-distributions ${\mathcal{D}^*_c}$ and ${\mathcal{D}^*}=\mathbb{I}{\mathcal{D}}$ are given by, $P_{\mathcal{D}^*_c}=e_3e_3^T$, and $P_{\mathcal{D}^*}=\left(I_{3\times 3}-P_{\mathcal{D}^*_c}\right)=-\widehat{e}_3^2$.
Denote by $\tau_\lambda$ the constraint moments that ensure these constraints. The Lagrange-d'Alembert principle implies that these constraint moments do no work. Thus$P_{\mathcal{D}^*}(\tau_{\lambda})\equiv0$ and hence the constraint moment must be of the form $\tau_{\lambda}=\left[0\:\: 0\:\: \tau_{{\lambda}_3}\right]^T$.
The external moments acting on the pendulum with respect to the pivot point of the pendulum is $\tau=\tau_u-Mgl\,e_3\times R^Te_3$. 
 Since $P_{\mathcal{D}^*_c}$ is constant, it follows that $\left(\nabla_{\Omega}P_{\mathcal{D}^*_c}\right)\left(\mathbb{I}\Omega\right)=e_3e_3^T\left(\mathbb{I}\Omega\times\Omega\right)$.
Thus from (\ref{eq:ConstrainetForce}) we have that the constraint moments are given by
\begin{align*}
\tau_{\lambda}=-e_3e_3^T\left(\mathbb{I}\Omega\times\Omega\right)=-
\begin{bmatrix}
0 \\ 0\\
(\mathbb{I}_1-\mathbb{I}_2)\Omega_1\Omega_2
\end{bmatrix},
\end{align*}
and from (\ref{eq:ConstrainedMech}) that the motion of the pendulum is described by
\begin{align*}
\dot{R}&=R\widehat{\Omega},\\
\mathbb{I}\nabla_{\Omega}\Omega&=-e_3e_3^T\left(\mathbb{I}\Omega\times\Omega\right)-Mgl\,{e}_3\times R^Te_3+\tau^u.
\end{align*}

We consider the control task of stabilizing the spherical pendulum in vertical upright position. Thus the output of interest is $y:SO(3) \mapsto \mathbb{S}^2$  that is explicitly given by $y(R)\triangleq Re_3$. Consider the Polar Morse function $V(y(R))=(1-e_3\cdot Re_3)$ on $\mathbb{S}^2$. The differential of this function is given by $dV=e_3\times Re_3$ and hence the Lie algebra version of the gradient of $V$ is given by $\mathbb{I}\eta=R^T\e_3\times e_3$. Note that $\eta_3\equiv 0$. Then since
\begin{align*}
\left(\nabla_{\Omega}P_{\mathcal{D}^*_c}\right)\left(\mathbb{I}\Omega_I\right)&=-\frac{1}{2} \left[\begin{array}{ccc}
0 \\
0 \\
\mathbb{I}_3(\Omega_1{\Omega_I}_2-\Omega_2{\Omega_I}_1) -(\mathbb{I}_1-\mathbb{I}_2)(\Omega_1{\Omega_I}_2+\Omega_2{\Omega_I}_1)
\end{array}\right].
\end{align*}
the PID controller (\ref{eq:ConstrianedI})--(\ref{eq:ConstrianedPID}) is explicitly given by
\begin{align}
\left[\begin{array}{cc}
{\tau_u}_1 \\ {\tau_u}_2
\end{array}\right]&=-k_p\left[\begin{array}{cc}\mathbb{I}_{1}\eta_1 \\ \mathbb{I}_{2}\eta_2
\end{array}\right]-k_d\left[\begin{array}{cc}\mathbb{I}_{1}\Omega_1 \\ \mathbb{I}_{2}\Omega_2
\end{array}\right]-k_I\left[\begin{array}{cc}\mathbb{I}_{1}\Omega_{I_1} \\ \mathbb{I}_{2}\Omega_{I_2}
\end{array}\right],  \label{eq:SP_PID} \\
\left[\begin{array}{cc}\mathbb{I}_{1}\dot{\Omega}_{I_1} \\ \mathbb{I}_{2}\dot{\Omega}_{I_2}
\end{array}\right]&=\left[\begin{array}{cc}\mathbb{I}_{1}\eta_1 \\ \mathbb{I}_{2}\eta_2
\end{array}\right]. \label{eq:SP_Integrator}
\end{align}

\subsubsection{Simulation Verification}
We simulate performance of the controller (\ref{eq:SP_PID})--(\ref{eq:SP_Integrator}) for the vertical stabilization of the spherical pendulum. The nominal parameters used for the controller were  $M=1\, kg$, $l=1\, m$, $g=1\, ms^{-2}$ and $\mathbb{I}=\mathrm{diag}\{1, 1, 1\}\, kgm^2$. To demonstrate the robustness of the controller the actual system parameters for the simulated system were chosen to be $50\%$ different from the nominal parameters used for the controller.

\begin{figure}[h!]
	\centering
	\begin{tabular}{ccc}
		\includegraphics[width=0.32\textwidth]{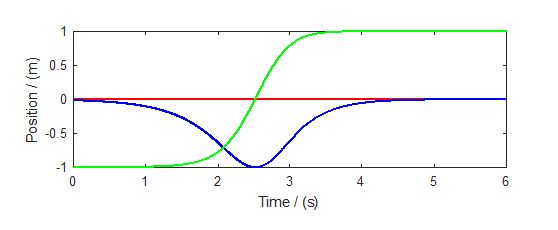} &
		\includegraphics[width=0.32\textwidth]{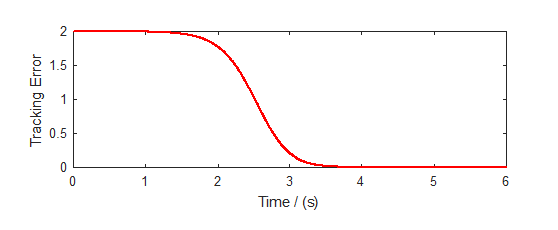} &
		\includegraphics[width=0.32\textwidth]{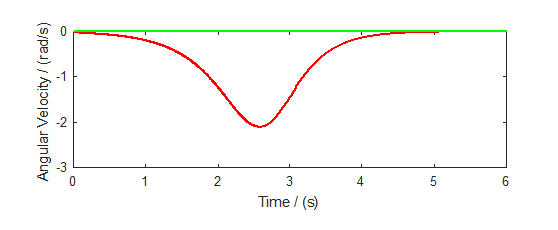}\\
		(a) & (b) & (c)
	\end{tabular}
	\caption{The figure shows the (a) position of the geometric center $o$, (b) tracking error $V(y(R))$ and (c) the body angular velocity $\Omega$ of the Spherical Pendulum with the PID controller (\ref{eq:SP_PID})--(\ref{eq:SP_Integrator}) in the presence of parameter uncertainties as large as $50\%$. \label{Fig:SP}}
\end{figure}

Figure-\ref{Fig:SP} presents simulation results for the position of the geometric center $o$, tracking error $V(y(R))$ and body angular velocity $\Omega$ of the Spherical Pendulum. Initial conditions used in the simulations were $R(0)=[1\:\:\:\: 0\:\:\:\: 0 ; 0\:\:\:\: \cos(179^0)\:\:\:\: -\sin(179^0) ; 0\:\:\:\: \sin(179^0)\:\:\:\: \cos(179^0)]$ and $\Omega(0)=[0\:\:\:\:0\:\:\:\:0]^T\,rads^{-1}$. These correspond to an initial condition very close to the vertically downward equilibrium of the system. The controller gains were chosen to be $k_p=16, k_d=8, k_I=1$.

\section{Conclusion}
This paper presents an intrinsic PID control framework for mechanical systems. It is motivated by the observation that a general unconstrained mechanical system evolving on a non-Euclidean manifold is essentially a double integrator system. Specifically we first show how to construct an intrinsic PID controller for fully actuated unconstrained mechanical system. On compact Lie groups the controller ensures almost-global locally exponential tracking of constant velocity references in the presence of unmodelled bounded constant disturbances and bounded parametric uncertainty. This controller is then extended to a class of interconnected underactuated mechanical systems. The stability properties of the resulting controller is first demonstrated for tracking of a hoop and then of a sphere rolling on a plane of unknown inclination.  
The third contribution is the extension to constrained mechanical systems on Lie groups. The stability properties of the controller is demonstrated using simulations that show the almost global and locally exponential stabilization of the vertically upright configuration of a spherical pendulum in the presence of bounded parametric uncertainty.


\begin{appendix}
\section{Proof of Theorem-\ref{theom:Theorem1}}
Consider the Morse function
\begin{align*}
W&=k_pV_s(V_y(g_s))+\frac{1}{2}\langle\langle{v_s}, {v_s}\rangle\rangle_s+\frac{\gamma}{2}\langle\langle {v_I}, {v_I}\rangle\rangle_s
+\alpha\langle\langle{\eta_e}, {v_s}\rangle\rangle_s+\beta \langle\langle {v_I}, {v_s}\rangle\rangle_s+\sigma \langle\langle {v_I}, {\eta_e}\rangle\rangle_s
+V_a(g_a)+\frac{1}{2}\langle\langle{v_a}, {v_a}\rangle\rangle_a
\end{align*}
Here $V_y : G_y \mapsto \mathbb{R}$ is a polar Morse function on $G_y$, and $V_a : G_a \mapsto \mathbb{R}$ is the polar Morse function corresponding to $\tau_V(g_a)$. Let ${\eta_e}$ be the gradient of $V_s\circ V_y$. When one chooses $\alpha=\frac{k_I}{k_d^2}$, $\beta=\frac{k_I}{k_d}$, $\sigma=2\kappa k_I$, $\gamma=\frac{k_I(k_I+k_pk_d)}{k_d^2}$,  $\mu=\mathrm{max}\,||\nabla\eta_e|| $ and $1/\mu<\kappa<2/\mu$. It can be shown that this function is bounded below by a quadratic function if (\ref{eq:kpCond}) is satisfied \cite{MaithripalaAutomatica}.

Differentiating $W$ along the dynamics of the closed loop system, we have
\begin{align*}
\dot{W}
&\leq-z_s^TQz_s-k_c\langle\langle {v_a},{v_a}\rangle\rangle_a+ \langle \tau_s+\Delta_s,{v_s}+\alpha{\eta_e}+\beta {v_I}\rangle
- \langle \epsilon_s(k_p{\eta_e}+k_d{v_s}+k_I{v_I}),{v_s}+\alpha{\eta_e}+\beta {v_I}\rangle
- k_c\langle \epsilon_a {v_a},{v_a}\rangle\\
&\:\:\:\:
-\langle k_cB^{-1}\mathbb{I}_a{v_a},{v_s}+\alpha{\eta_e}+\beta {v_I}\rangle
-\langle B\mathbb{I}_s(k_p{\eta_e}+k_d{v_s}+k_I{v_I}),{v_a}\rangle
+\langle  \tau_{a}+\Delta_a, {v_a}\rangle\\
& \leq -(\lambda_{\mathrm{min}}(Q)-g_0g_1||\epsilon_{s}||)||z_s||^2
+g_1||\Delta_s||(||{v_I}|| + ||{\eta_e}|| + ||{v_s}||)
-((k_c\lambda_{\mathrm{min}}(\mathbb{I}_a)-k_c||\epsilon_{a}||)||{v_a}||-||\Delta_a||)||{v_a}||\\
&\:\:\:\: +g_1g_3||{v_a}||||z_s||^2+g_4||z_s||||{v_a}||^2
+g_5||\epsilon_A||||{v_a}||(||{v_I}|| + ||{\eta_e}|| + ||{v_s}||)
\end{align*}
where
 \begin{align*}
 Q&=\begin{bmatrix}
 \frac{k_I^2}{k_d} & 0 & -\delta k_I\\
 0 &   \frac{k_I}{k_d^2}(k_p-2\kappa k_d^2) & 0 \\
-\delta k_I & 0 & k_d- \frac{\mu k_I}{k_d^2}
 \end{bmatrix},\:\:\:\:\:\:\:
 z_s=\begin{bmatrix}
||{v_I}|| & ||{\eta_e}|| & ||{v_s}||\end{bmatrix},
 \end{align*}
$g_0=\mathrm{max}\{k_p,k_d,k_I\}$, $g_1=\mathrm{max}\{1,k_I/k_d,k_I/k_d^2\}$, 
$g_3=\max \frac{\tau_s({v_s},{v_a})}{||{v_s}||||{v_a}||}$, $g_4=\max \frac{\tau_a({v_s},{v_a})}{||{v_s}||||{v_a}||}$, $\epsilon_{s}$  and $\epsilon_{a}$ are operators that depend on the magnitude of the parameter uncertainties. Here $\lambda_{\mathrm{min}}(A)$ is the smallest eigenvalue of the matrix $A$. It can be shown that $Q$ is positive definite if the gains are chosen such that (\ref{eq:kICond})--(\ref{eq:kpCond}) are satisfied.

Given any compact subset $\mathcal{X}$, and a small $\epsilon>0$, let $\mathcal{W}_u$ be the smallest level set of $W$ containing  
$\mathcal{X}$ and let $k_{\mathcal{X}}>0$ be the smallest value such that $||{v_I}|| + ||{\eta_e}|| + ||{v_s}||+||{v_a}||<k_{\mathcal{X}}$ for all $(y,{v_s},{v_I},g_a,{v_a})\in \mathcal{W}_u$ and let $\mathcal{W}_l$ be the smallest level set of $W$ containing the set where $||z||<\epsilon$. 
On $\mathcal{W}_u$ we find that
\begin{align*}
\dot{W}
& \leq -(\lambda_{\mathrm{min}}(Q)-g_0g_1||\epsilon_{s}||)||z_s||^2
+g_1||\Delta_s||(||{v_I}|| + ||{\eta_e}|| + ||{v_s}||)
-(k_c\lambda_{\mathrm{min}}(\mathbb{I}_a)-k_c||\epsilon_{a}||)||{v_a}||^2+||\Delta_a||||{v_a}||\\&\:\:\:\: 
+k_{\mathcal{X}}g_1g_3||z_s||^2+k_{\mathcal{X}}g_4||{v_a}||^2
+g_5||\epsilon_A||||{v_a}||(||{v_I}|| + ||{\eta_e}|| + ||{v_s}||)
\\& \leq -(\lambda_{\mathrm{min}}(\tilde{Q})-g_7||\epsilon_A||)||z||^2
+g_6(||{v_I}|| + ||{\eta_e}|| + ||{v_s}||+||{v_a}||)
\end{align*}
where 
\begin{align*}
g_6&\triangleq\max \{g_1||\Delta_s||,||\Delta_a||\},\\
\lambda_{\mathrm{min}}(\tilde{Q})&\triangleq\min \{
(\lambda_{\mathrm{min}}(Q)-g_0g_1||\epsilon_{s}||-k_{\mathcal{X}}g_3),\\
&\:\:\:\: (k_c\lambda_{\mathrm{min}}(\mathbb{I}_a)-k_c||\epsilon_a||-k_{\mathcal{X}}g_4)\},
\end{align*}
and $z=[||{v_I}|| \:\:\:\: ||{\eta_e}|| \:\:\:\: ||{v_s}||\:\:\:\: ||{v_a}||]^T$.

Given any $k_{\mathcal{X}}$ and $\epsilon >0$ let the gains be picked sufficiently large enough so that
\begin{align*}
(\lambda_{\mathrm{min}}(\tilde{Q})-g_7||\epsilon_A||)&>0,\\
\frac{g_6}{(\lambda_{\mathrm{min}}(\tilde{Q})-g_7||\epsilon_A||)}&<{\epsilon}^2.
\end{align*}
Then we see that $W$ is strictly decreasing in $\mathcal{W}_u/\mathcal{W}_l$. Thus the Lasalle's invariance principle implies that the trajectories of (\ref{eq:System_s}), (\ref{eq:System_a}), and (\ref{eq:GeneralI}) converge to the largest invariant set contained in the set where $\dot{W}\equiv 0$ contained in $\mathcal{W}_u/\mathcal{W}_l$. For mechanical systems with constant velocity references and constant unknown disturbances these invariant sets are exactly the equilibria of (\ref{eq:System_s}) -- (\ref{eq:System_a}). Thus proving that $\lim_{t\to \infty} y(t)=0$ and $\lim_{t\to \infty}({v_s}(t),{v_a}(t))=(0,0)$ semi-globally in the presence of bounded parametric uncertainty and bounded constant disturbances for all initial conditions in $\mathcal{X}$ other than the unstable equilibria and their stable manifolds. The local exponential stability follows from the fact that $W$ is quadratically bounded from below.


\end{appendix}

\bibliographystyle{elsarticle-num}
\bibliography{DHS_Maithripala}

\end{document}